\documentclass[twocolumn]{svjour3}          

\usepackage{amsmath,amssymb,amsfonts,xcolor}
\usepackage{multirow}
\usepackage{booktabs}
\usepackage{ifthen}
\usepackage{xspace}
\usepackage{subfigure}
\usepackage[pdftex]{graphicx}

\usepackage{tikz}
\usetikzlibrary{patterns}
\usetikzlibrary{matrix}
\usetikzlibrary{positioning}

\newtheorem{algorithm}{\textbf{\textup{Algorithm}}}

\let\thetaTMP\theta
\let\theta\vartheta
\let\vartheta\thetaTMP
\newcommand{\Alg}{\mathcal A}

\newcommand{\PDx}{\mathcal{PD}_u}
\newcommand{\PDy}{\mathcal{PD}_p}
\newcommand{\enquote}[1]{``#1''}
\newcommand{\dom}{\operatorname{dom}}
\newcommand{\cl}[1]{\operatorname{cl}(#1)}
\renewcommand{\int}[1]{\operatorname{int}(#1)}
\newcommand{\bdry}[1]{\operatorname{bdry}(#1)}
\newcommand{\ind}[1]{\ensuremath{\delta_{#1}}}

\newcommand{\R}{\mathbb R}

\newcommand{\N}{\mathbb N}
\newcommand{\eR}{\overline{\R}}
\newcommand{\st}{s.t.\xspace}
\newcommand{\map}[3]{#1\colon #2 \to #3}
\newcommand{\norm}[2][]{\Vert #2 \Vert_{#1}}
\newcommand{\scal}[2]{\left\langle #1,#2 \right\rangle}
\newcommand{\seq}[2][]{(#2)_{#1}}
\newcommand{\set}[2][]{\ifthenelse{\equal{#1}{\empty}}{\lbrace #2 \rbrace}{\lbrace #2\vert\, #1\rbrace}}
\newcommand{\abs}[1]{\vert #1 \vert}

\newcommand{\proj}{\operatorname{proj}}
\newcommand{\prox}[1][]{\operatorname{prox}^{#1}}
\newcommand{\ncone}[1]{\mathrm{N}_{#1}}

\newcommand{\quadand}{\quad\text{and}\quad}

\newcommand{\Simplex}{\Delta}
\newcommand{\cost}{\mathfrak{c}}
\newcommand{\cnn}{\mathfrak{f}}
\newcommand{\nx}{{N_x}} 
\newcommand{\ny}{{N_y}} 
\newcommand{\nz}{{N_l}} 
\newcommand{\dimOpt}{N} 
\newcommand{\dimP}{P} 
\newcommand{\aitvar}{n} 
\newcommand{\ait}[1]{{(#1)}} 
\newcommand{\nc}{{N_c}} 

\newcommand{\px}{i}
\newcommand{\ppi}{\mathbf{i}} 
\newcommand{\ppd}{\mathbf{j}} 
\newcommand{\py}{j}
\newcommand{\pz}{k}

\newcommand{\optu}{u} 
\newcommand{\optp}{p} 
\newcommand{\opK}{\nabla} 
\newcommand{\opKnz}{\mathcal D} 
\newcommand{\WW}{\mathcal W} %

\newcommand{\pt}{t} 
\newcommand{\nt}{{N_T}} 
\newcommand{\img}{\mathfrak I} 
\newcommand{\gt}{\mathfrak g} 

\newcommand{\Loss}{\mathcal L} 
\newcommand{\loss}{\ell} 
\newcommand{\Hess}{H} 

\newcommand{\X}{{X}}
\newcommand{\Y}{{Y}}

\newcommand{\titemi}{(i)~}%
\newcommand{\titemii}{(ii)~}%
\newcommand{\titemiii}{(iii)~}%

\definecolor{dSkyBlue}{rgb}{.11,.2,.62}
\definecolor{bSkyBlue}{rgb}{.13,.73,.92}

\usepackage{cancel}
\usepackage[normalem]{ulem}

\newif\ifShowChanges
\ShowChangestrue
\ifShowChanges

\else

\fi

\def\mygraphicspath{./figures/}

\smartqed  
\usepackage{graphicx}
\begin{document}

\title{Techniques for Gradient Based\\ Bilevel Optimization with\\ Nonsmooth Lower Level Problems}
\titlerunning{Bilevel Optimization with Nonsmooth Lower Level Problems} 

\author{Peter Ochs \and
        Ren\'e Ranftl \and \\
        Thomas Brox \and
        Thomas Pock
}


\institute{P. Ochs \at
              Mathematical Image Analysis Group\\
              University of Saarland, Germany \\
              \email{ochs@mia.uni-saarland.de}
           \and
           R. Ranftl \at
            Visual Computing Lab\\
            Intel Labs, Santa Clara, CA, United States\\
            \email{rene.ranftl@intel.com}        
           \and
           T. Brox and P. Ochs\at
              Computer Vision Group\\
              University of Freiburg, Germany \\
              \email{\{ochs,brox\}@cs.uni-freiburg.de}
           \and
           T. Pock \at
              Institute for Computer Graphics and Vision\\
              Graz University of Technology, Austria\\      
              and\\
              Digital Safety \& Security Department\\
              AIT Austrian Institute of Technology GmbH \\
              1220 Vienna, Austria\\
              \email{pock@icg.tugraz.at}
}

\date{Received: date / Accepted: date}

\maketitle

\begin{abstract}
  We propose techniques for approximating bilevel optimization problems with non-smooth  lower level problems that can have a non-unique solution. To this end, we substitute the expression of a minimizer of the lower level minimization problem with an iterative algorithm that is guaranteed to converge to a minimizer of the problem. Using suitable non-linear proximal distance functions, the update mappings of such an iterative algorithm can be differentiable, notwithstanding the fact that the minimization problem is non-smooth. 
\end{abstract}

\newpage

\section{Introduction}

We consider numerical methods for solving bilevel optimization problems of the form
\begin{equation} \label{eq:bilevel-general-intro}
\begin{split}
  &\ \min_\theta\ \Loss(x^*(\theta), \theta) \\
  &\ \st\ x^*(\theta) \in \arg\min_{x\in\R^\dimOpt} E(x,\theta) \,,
\end{split}
\end{equation}
where $\Loss$ (denoted \emph{loss function}) is a function penalizing the differences between the output of the lower level problem $x^*(\theta)$ and some given ground truth data. In addition, $\Loss$ can also contain a regularizer on the parameter vector $\theta$, e.g. a sparsity prior. The mapping $x^*(\theta)$ is the solution of an optimization problem (parametrized by $\theta$) that solves a specific task, e.g. multi-label segmentation.

In the general bilevel literature, \eqref{eq:bilevel-general-intro} is often presented as a leader--follower problem. The leader (upper level problem) tries to optimize the next move (minimization of the upper level problem) under consideration of the move of an opponent, the follower. Given some information $\theta$ to the follower, the leader tries to anticipate the follower's next move (minimization of the lower level problem).

In this paper, we focus on a class of problems that allows for non-smooth convex functions $x\mapsto E(x,\theta)$ in the lower level problem, e.g. sparse models based on the $\ell_1$-norm.  Such models have become very popular in the computer vision, image processing and machine learning communities since they are robust with respect to noise and outliers in the input data. 

Due to the possibly high dimensionality of the parameter vector, we pursue the minimization of the bilevel problem \eqref{eq:bilevel-general-intro} using gradient based methods. Hence, a descent direction of $\Loss$ with respect to $\theta$ must be determined. Its estimation involves the Jacobian of the solution map $x^*(\theta)$ with respect to the parameter vector $\theta$, which causes three kinds of problems:
\begin{itemize}
\item[\titemi] The solution mapping $x^*(\theta)$ is only defined implicitly (as a minimizer of the lower level problem).
\item[\titemii] The lower level's solution is not unique.\footnote{Note that the bilevel problem as in \eqref{eq:bilevel-general-intro} is not well-defined in this case. We discuss some details in Section~\ref{sec:bilevel-problem}.}
\item[\titemiii] The lower level problem is non-smooth.
\end{itemize}

\titemi %
A reduction to a single level problem by explicitly solving the lower level problem is not always possible. Nevertheless, if the lower level problem is sufficiently smooth, sometimes, it can be replaced by its optimality condition, and the implicit function theorem (cf. Section~\ref{sec:derivative-impl-fun}) provides an explicit formula for the derivative of the solution map. This approach does not work for non-smooth lower level problems.

\titemii Consider the example 
\begin{equation} \label{eq:bilevel-example-nonunique-A}
\begin{split}
  &\ \min_{\theta\in \R} \ (x^*(\theta)-1)^2 \\
  &\ \st\ x^*(\theta) \in \arg\min_{x\in [0,1]} \theta x \,,
\end{split}
\end{equation}
which reduces to minimization of a step function
\[
  \min_{\theta\in\R}\, \Loss(x^*(\theta))\,, \quad \Loss(x^*(\theta)) = 
    \begin{cases} 1,& \text{if } \theta > 0\,; \\ 
                  0,& \text{if } \theta < 0\,; \\ 
              [0,1],& \text{if } \theta = 0 \,. \end{cases}
\]

A gradient based method will get stuck almost everywhere, as the derivative vanishes for all $\theta \neq 0$. Similar situations arise for robust models in the lower level problem. By definition the solution is not affected by small perturbations of the input data (or the parameter $\theta$). For instance in the multi-label segmentation problem, small changes in the pixel likelihoods do not change the segmentation result; the energy landscape of the loss function will have the form of a high dimensional step function. 

\titemiii Due to the non-smoothness of the lower level problem, standard calculus cannot be applied. In variational (non-smooth) analysis, there are many generalizations of derivatives, such as the convex subdifferential, the Fr\'echet subdifferential, or the limiting subdifferential. However they are often set-valued and generalizations of the chain rule and rely on constraint qualifications that are sometimes quite restrictive and often hard to verify. \\

In the conference version of this paper \cite{ORBP15}, we introduced an approach to overcome the smoothness restriction in some cases of practical interest. The idea is to replace the lower level problem by an iterative algorithm that is guaranteed to converge to a solution of the problem. If the update mapping of the algorithm is a smooth function, the chain rule can be applied to the composition of these update mappings recursively and the exact derivatives with respect to the parameter vector $\theta$ can be computed. Algorithms based on Bregman distances are key for this development. The number of iterations of the iterative algorithm steers the approximation quality of the lower level problem.

The iterative algorithm that replaces the lower level is stopped after a small number of iterations. However, once the algorithm and the number of iterations are fixed, the resulting bilevel optimization problem seeks for optimal $\theta$ for exactly this algorithm and this (fixed) number of iterations. Numerically, the derivative that is involved in gradient based minimization is exact: the number of chain rule recursions is finite.
This is in contrast to an approach based on the optimality condition of a smooth approximation of the lower level problem. In this case, the descent direction is based on the derivative of the optimality condition evaluated at the minimizer of the lower level problem, which is only determined approximately.

Beyond the analysis of the conference paper, we discuss approximations to the derivative evaluation that reduce the memory requirements and the computational cost significantly. We extend the class of problems that can be used in our framework and give some more details about the general implementation of our approach. Moreover, we consider the limiting case, i.e., the fixed point equation of an iterative algorithm in the lower level problem.

We point out several applications of our approach and evaluate it for a multi-label segmentation problem  coupled with a convolutional neural network.

\section{Related Work}

We propose a simple approximation of the lower level problem that naturally addresses \emph{non-smoothness} and \emph{non-uniqueness}. 

For a non-unique solution map (a set-valued mapping) of the lower level problem \eqref{eq:bilevel-general-intro} is not even well-defined (cf. Remark~\ref{rem:well-posed-not-unique}). \cite{Dempe15} describes three possible options to cope with this problem. The \emph{optimistic bilevel optimization problem} assumes a cooperative strategy of leader and follower, i.e., in case of multiple solutions the follower tries to minimize the upper level objective. The \emph{pessimistic bilevel problem} is the other extreme, where the leader tries to bound the damage that the follower could cause by its move. The \emph{selection function approach} assumes that the leader can always predict the followers choice. Of course, these three approaches are the same for lower level problems with a unique output.

Our approach does not fall into any of the three cases, however the selection function approach is the closest. The difference is that our approximation changes the output also at (originally) unique points. Our solution strategy reduces the solution map to be single-valued, similar to the approaches mentioned above.

Dempe et al. \cite{Dempe15} classifies the following  optimality conditions\footnote{The classification in \cite{Dempe15} applies to the optimistic bilevel problem.}. The \emph{primal Karush--Kuhn--Tucker (KKT) transformation} replaces the lower level problem by the necessary and sufficient optimality condition for a convex function. The equivalence to the original problem is shown in \cite{DZ10}. The \emph{classical KKT transformation} substitutes the lower level problem with the classical KKT conditions. Due to the  extra variable, the problems are not fully equivalent anymore (see \cite{Dempe15}). This approach, which leads to a non-smooth mathematical problem with complementary constraints (MPEC), is the most frequently used one. The third approach is the \emph{optimal value transform}, which introduces a constraint that bounds the lower level objective by the optimal value function.

Our approach is---in the limit---motivated by the first class of the primal KKT transformation. We consider the fixed point equation of an algorithm, which represents the optimality condition without introducing additional variables, and approximate this situation with finitely many iterations of the algorithm. 

We focus on gradient based methods, such as gradient descent, L-BFGS \cite{BFGS89}, non-linear conjugate gradient \cite{FR64,Al-Baali85}, Heavy-ball method \cite{ZK93}, iPiano \cite{OCBP14}, and others, for solving the bilevel optimization problem. 
In particular, this paper focuses on the estimation of descent directions. As one option, the gradient can be approximated numerically with finite differences such as in \cite{Domke10}. We rather pursue what is known as algorithmic/automatic differentiation. It is based on the idea to decompose the derivative evaluation into small parts by means of a chain rule, where the analytic derivative of each part is known. A whole branch of research deals with this technique \cite{GW08}. Obviously, the idea to differentiate an algorithm in the lower level problem is not new \cite{Tappen07,Domke12}. The difference is that our algorithm has a smooth update mapping while actually minimizing a non-smooth objective function. Another idea to approach a non-smooth problem with an iterative algorithm is presented in \cite{DVFP14}, where a chain rule for weak derivatives is used (cf. Section~\ref{sec:gen-derivative-alg}).


The special case of a lower level problem that depends linearly on the parameters is treated by structured output support vector machines~\cite{Tsochantaridis05}. The linear structure of the lower level problem allows the construction of an upper bound of the upper level objective function, which needs to be minimized. In general, this approach is only an approximation to the bilevel problem in \eqref{eq:bilevel-general-intro}, which can be solved using subgradient descent.

There are several practical examples of bilevel optimization in the computer vision and machine learning. Bilevel optimization was considered for task specific sparse analysis prior learning \cite{PF11} and applied to signal restoration. In \cite{KP13,CPRB13,CRP14} a bilevel approach was used to learn a model of natural image statistics, which was then applied to various image restoration tasks. A variational formulation for learning a good noise model was addressed in \cite{RS13} in a PDE-constrained optimization framework, with some follow-up works \cite{CRS14,RSV15,CRSV15}. In machine learning bilevel optimization was used to train a SVM \cite{BKHP08} and other techniques \cite{Moore10}. Recently, it was used for the end-to-end training of a Convolutional Neural Network (CNN) and a graphical model for binary image segmentation \cite{RP14} (cf. Section~\ref{sec:potts-seg-CNN}).


Finally, we  refer to \cite{Dempe03} for an annotated bibliography with many references regarding the theoretical and practical development in bilevel optimization.

\section*{Preliminaries} 

We work in a Euclidean vector space $\R^\dimOpt$ of dimension $\dimOpt$ equipped with the standard Euclidean norm $\norm\cdot := \sqrt{\scal\cdot\cdot}$ that is induced by the standard inner product. We use the notation $\eR:=\R\cup\set{\infty}$ to denote the extended real numbers.

\renewcommand{\bullet}{*}
We use the notation $[x \bullet a]$ for $x,a\in \R^\dimOpt$ to denote the set $\set{x\in\R^\dimOpt\vert\, \forall i\colon\, x_i \bullet a_i}$, where $\bullet\in\set{<,\leq, =, \geq, >}$ is a binary relation on $\R\times\R$. For example $[x\geq 0]$ denotes the non-negative orthant in $\R^\dimOpt$.

\section{The Bilevel Problem} \label{sec:bilevel-problem}

We consider bilevel optimization problems of the form:
\begin{equation} \label{eq:bilevel-general}
\begin{split}
  &\ \min_{\theta\in\R^\dimP}\ \Loss(x^*(\theta), \theta) + \loss(\theta)\\
  &\ \st\ x^*(\theta) \in \arg\min_{x\in\R^\dimOpt} E(x,\theta)
\end{split}
\end{equation}
The function $\map \loss {\R^\dimP}{\eR}$ is assumed to be proper, lower semi-continuous, convex, and \enquote{prox-friendly}\footnote{The associated proximity operator has a closed-form solution or the solution may be determined efficiently numerically.} and the function $\map \Loss{\R^\dimOpt\times\R^\dimP}{\R}$ to be continuously differentiable on $\dom\loss$. The optimization variable is the (parameter) vector $\theta\in\R^\dimP$. It appears implicitly and explicit in the upper level problem. It is implicit via the solution mapping $x^*(\theta)\in\R^\dimOpt$ of the lower level problem and explicit in $\loss$ and in the second argument of $\Loss$. The lower level is a minimization problem in the first variable of a proper, lower semi-continuous function $\map{E}{\R^\dimOpt\times\R^\dimP}{\eR}$. For each $\theta\in\R^\dimP$ the objective function (energy) $x \mapsto E(x,\theta)$ is assumed to be convex.\\
Note that our formulation includes constrained optimization problems in the upper and lower level problem. The functions $\loss$ and $E$ are defined as extended-valued (real) functions. Of course, in order to handle the constraints efficiently in the algorithm, the constraint sets should not be too complicated.

In order to solve the optimization problem in \eqref{eq:bilevel-general}, we can apply iPiano \cite{OCBP14}, a gradient-based algorithm that can handle the non-smooth part $\loss(\theta)$. The extension of iPiano in \cite[Chapter 6]{Ochs15} allows for a prox-bounded (non-convex, non-smooth) function $\loss(\theta)$. Informally, the update step of this algorithm (for the parameter vector $\theta$) reads
\begin{multline} \label{eq:ipiano-update}
  \theta^{k+1} \in \prox_{\alpha_k \loss}\left(\theta^{k} - \alpha_k \nabla_\theta \Loss(x^*(\theta^k), \theta^k) \right.\\
    \left.+ \beta_k(\theta^{k} - \theta^{k-1})\right) \,,
\end{multline}
where $\prox_{\alpha_k \loss}$ denotes the proximity operator of the function $\loss$, and $\alpha_k$ is a step-size parameter and $\beta_k$ steers the so-called inertial effect of the algorithm (usually $\beta_k\in[0,1]$). For details about $\alpha_k$ and $\beta_k$, we refer to \cite{OCBP14,Ochs15}, where convergence to a stationary point (a zero in the limiting subdifferential) is proved under mild assumptions. We could also apply proximal gradient descent (forward--backward splitting) \cite{ABS13} ($\beta_k=0$). In our experiments, iPiano was usually faster and less sensitive to local optima, however. If the non-smooth term is not present, several gradient based solvers can be used \cite{BFGS89,ZK93,FR64,Al-Baali85}. 

The structure of the update step in \eqref{eq:ipiano-update} points out that the main aspect in applying such gradient-based algorithm\xspace is the evaluation of the gradient $\nabla_\theta \Loss(x^*(\theta), \theta)$. The remainder of this paper deals with exactly this problem: compute $\nabla_\theta \Loss(x^*(\theta), \theta)$ with a solution mapping $\theta\mapsto x^*(\theta)$ of a possibly non-smooth objective function in the lower level. Note that the following approximations naturally yield or require a unique solution of the lower level problem.

\begin{remark} \label{rem:well-posed-not-unique}
The formulation \eqref{eq:bilevel-general} of a bilevel optimization problem only makes sense when $\arg\min_{x\in\R^\dimOpt} E(x,\theta)$ yields a unique minimizer. In that case, optimality of the bilevel problem can be derived from standard optimality conditions in non-linear programming. If the lower level problem does not provide a unique solution, the loss function $\Loss$ must be defined on the power set of $\R^\dimOpt$ and a different notion of optimality must be introduced. Since, this results in problems beyond the scope of this paper, we refer to \cite{Dempe15}. A common circumvention is to consider the corresponding optimistic bilevel problem.
\end{remark}

\section{Computing descent directions} \label{sec:descent-directions}

For a given parameter value $\theta\in\R^\dimP$, we would like to compute a descent direction of $\Loss$ in \eqref{eq:bilevel-general} with respect to $\theta$ to find a numerical solution using some gradient based method. Obviously, we need the derivative of the solution map $x^*(\theta)$ with respect to $\theta$. In the following, we present strategies to approximate the (possibly non-smooth) lower level problem and to compute a descent direction.

\subsection{Derivative of a smoothed lower level problem} \label{sec:derivative-impl-fun}

If  the objective function of the lower level problem of \eqref{eq:bilevel-general} can be approximated well with a twice continuously differentiable function (again denoted $E$), we can make use of the implicit function theorem to find the derivative of the solution map with respect to $\theta$. The optimality condition of the lower level problem is $\nabla_x E(x,\theta) = 0$, which under some conditions implicitly defines a function $x^*(\theta)$.  As we assume that the problem $\min_x E(x,\theta)$ has a solution, there is $(x^*,\theta)$ such that $\nabla_x E(x^*,\theta) = 0$. Then, under the conditions that $\nabla_x E(x^*,\theta)$ is continuously differentiable and $({\partial (\nabla_x E)}/{\partial x})(x^*,\theta)$ is invertible, there exists an explicit function $x^*(\theta)$ defined on a (open) neighborhood of $x^*$. Moreover, the function $x^*(\theta)$ is continuously differentiable at $\theta$ and it holds that
\[
    \frac{\partial x^*}{\partial \theta} (\theta) =\! \left( - \frac{\partial (\nabla_x E)}{\partial x}(x^*(\theta),\theta) \right)^{-1}\! \frac{\partial (\nabla_x E)}{\partial \theta} (x^*(\theta),\theta) \,.
\]
Using the Hessian $\Hess_E(x^*( \theta), \theta) := \frac{\partial^2 E}{\partial x^2}(x^*( \theta), \theta)$ yields
\begin{equation} \label{eq:implicit-diff-param}
    \frac{\partial x^*}{\partial \theta} ( \theta) = - (\Hess_E(x^*( \theta), \theta))^{-1} \frac{\partial^2 E }{\partial \theta\partial x} (x^*( \theta), \theta) \,.
\end{equation}
The requirement for using \eqref{eq:implicit-diff-param} from the implicit function theorem is the continuous differentiability of ${\partial E}/{\partial x}$ and the invertibility of $\Hess_E$. Application of the chain rule yields the total derivative of the loss function $\Loss$ of \eqref{eq:bilevel-general} w.r.t. $\theta$
\begin{equation} \label{eq:full-derivative-impl-diff-param}
  \frac{d\Loss}{d\theta} = - \Bigg[ \frac{\partial \Loss }{\partial x} \Hess_E^{-1}\Bigg] \frac{\partial^2 E }{\partial \theta\partial x} + \frac{\partial\Loss}{\partial\theta}\,,
\end{equation}
where the function evaluation at $(x^*(\theta),\theta)$ is dropped for brevity. A clever way of setting parentheses, as it is indicated by the squared brackets, avoids explicit inversion of the Hessian matrix. However, for large problems iterative solvers are required.

\subsection{Derivative of iterative algorithms} \label{sec:derivative-iterative-alg}

We can replace the minimization problem in the lower level of \eqref{eq:bilevel-general} by an algorithm that solves this problem, i.e., the lower level problem is replaced by an equality constraint. This approach shows three advantages: \titemi After approximating the lower level of \eqref{eq:bilevel-general} by an algorithm, the approach is exact; \titemii the update step of the algorithm can be smooth without the lower level problem to be smooth; \titemiii the output is always unique (for a fixed initialization), which circumvents the critical issue of a non-unique lower level solution.

Let $\map{\Alg\text{ and }\Alg^\ait{\aitvar}}{X\times \R^\dimP}{X}$ describe one or $\aitvar$ iterations, respectively, of algorithm $\Alg$ for minimizing $E$ in \eqref{eq:bilevel-general}. For simplicity, we assume that the feasible set mapping $\theta\mapsto \set{x\in\R^\dimOpt\vert\, (x,\theta)\in \dom\Loss}$ is constant\footnote{More generally, the concept of outer semi-continuity of the feasible set mapping is needed, otherwise a gradient based method could converge to a non-feasible point.}, i.e., the same $X$ is assigned to all $\theta\in\R^\dimP$. Note that $X=\R^\dimOpt$ is permitted. 

For  fixed $\aitvar\in\N$, we replace \eqref{eq:bilevel-general} by 
\begin{equation}\label{eq:bilevel-abstract-alg}
\begin{split}
  &\ \min_\theta\ \Loss(x^*(\theta), \theta)  + \loss(\theta) \\
  &\ \st\ x^*(\theta) = \Alg^\ait{\aitvar+1} (x^\ait0, \theta) \,,
\end{split}
\end{equation}
where $x^\ait0$ is some initialization of the algorithm. The solution map of the lower level problem $x^*(\theta)$ is the output of the algorithm $\Alg$ after $\aitvar+1$ iterations. If we write down one iteration of the algorithm, i.e., $x^\ait{\aitvar+1}(\theta) = \Alg(x^\ait\aitvar(\theta),\theta)$, we have to assume that $x^\ait\aitvar$ depends on the choice of $\theta$. However, this dependency can be dropped for the first iterate, which emerges from the initialization. \\

A suitable algorithm has the properties that $x^\ait\aitvar(\theta)$ converges pointwise (i.e. for each $\theta$) to a solution of the lower level problem as $\aitvar$ goes to infinity and $E(x^\ait\aitvar,\theta) = E(\Alg^\ait\aitvar (x^\ait0, \theta),\theta) \to \min_x E(x,\theta)$ for $n\to\infty$. Note that for Bregman proximity functions in  algorithm $\Alg$, the solution for $n\to\infty$ could lie on $\bdry X$, despite $x^\ait\aitvar\in\int X$ for all $\aitvar$. However, this matters only for an asymptotic analysis.\\

 If $\Alg$ is (totally) differentiable with respect to $\theta$, then, by the standard chain rule, $\Alg^\ait\aitvar$ is differentiable with respect to $\theta$ as well. This way, we obtain a totally differentiable approximation to the lower level problem of \eqref{eq:bilevel-general}, where the approximation quality can simply be controlled by the number of iterations. 
For so-called descent algorithms, it holds that
\[
  E(x^\ait{\aitvar+1},\theta) - \min_x E(x,\theta) \leq E(x^\ait{\aitvar},\theta) - \min_x E(x,\theta) \,.
\]
A large number of iterations usually approximates the minimum of $E$ better than a small number of iterations. 

Nevertheless, also a small number of iterations is interesting for our approach. Once a certain number of iterations is fixed, the bilevel optimization problem seeks for an optimal performance with exactly this chosen number of iterations. Solving the bilevel optimization problem accurately with a small number of iterations $\aitvar$ of the lower level algorithm can result in a better performance than a poorly solved  bilevel problem with a large number of iterations in the lower level. 

Our approach is well suited for minimizing the bilevel problem using gradient based methods. The differentiation of $\Loss$ with respect to $\theta$ in \eqref{eq:bilevel-abstract-alg} is exact; one an algorithm is selected no additional approximation is required for computing the derivatives. In contrast, the smoothing approach from Section~\ref{sec:derivative-impl-fun} requires the minimization of a smooth objective function, the solution of which can be found only approximatively. Therefore, the descent direction, which is based on the optimality condition, is always erroneous.

The \enquote{smoothing parameter} in our approach is the number of iterations of the algorithm that replaces the lower level problem. Since the algorithm's update mapping is assumed to be smooth, in particular, locally Lipschitz continuous, which formally means
\[
  \norm{\Alg(x, \theta) - \Alg(y, \theta)} \leq \operatorname{const.} \norm{x-y}
\]
holds in a neighborhood of the initial point, the variation of the output after one iterations is limited.  Therefore, intuitively, for a large number of iterations $\aitvar$, less smoothness of $\Alg^\ait\aitvar$ can be expected. \\

In order to obtain the derivative of the lower level problem of \eqref{eq:bilevel-abstract-alg}, there are two prominent concepts: forward mode and backward mode. 
For any vector $\xi \in\R^\dimOpt$, the \emph{forward mode} corresponds to evaluating the derivative as 
\begin{multline} \label{eq:alg-iter-fwd}
  \xi^\top \frac{dx^\ait{\aitvar+1}}{d\theta}(\theta) = \\
  \xi^\top \left[ \frac{\partial\Alg}{\partial x}(x^\ait\aitvar,\theta) \frac{dx^\ait{\aitvar}}{d\theta}(\theta) \right]+ \xi^\top \frac{\partial \Alg}{\partial \theta} (x^\ait\aitvar,\theta)\,,
\end{multline}
whereas the \emph{backward mode/reverse mode} evaluates the derivative as 
\begin{multline} \label{eq:alg-iter-reverse}
  \Bigg(\frac{dx^\ait{\aitvar+1}}{d\theta}(\theta)\Bigg)^\top \xi \\= \Bigg( \frac{dx^\ait{\aitvar}}{d\theta}(\theta)\Bigg)^\top \left[ \Bigg(\frac{\partial\Alg}{\partial x}(x^\ait\aitvar,\theta)\Bigg)^\top \xi \right]\\ + \left( \Bigg(\frac{\partial \Alg}{\partial \theta} (x^\ait\aitvar,\theta)\Bigg)^\top \xi \right) \,,
\end{multline}
where the squared brackets symbolize the different orders of evaluating the terms. In both approaches, replacing and evaluating the term ${dx^\ait{\aitvar}}/{d\theta}$ using the preceding iterate $\ait{\aitvar-1}$ is done in the respective order.

Mathematically both concepts result in the same solution. However, numerically the approaches are very different. The reverse mode is usually more efficient when the optimization variable $\theta$ is high dimensional (i.e., $\dimP$ is large) and the range of the objective function $\Loss$ is low dimensional---it is always $1$ in our setting. This corresponds to $\xi$ being a column vector instead of a derivative matrix. 
The forward mode is often easier to implement, since it is executed in the same order as the optimization algorithm itself and can be computed online, i.e., during the iteration of the algorithm. As a downside, each partial derivative must be initialized and propagated through the iterations. Therefore, the memory requirement is vastly increasing with the dimension $\dimP$. We focus on the reverse mode for evaluating the derivatives, due to its computationally more appealing nature.

The backward mode is executed in the reverse order of the iterations of the algorithm and needs the optimum $x^*$, which is $x^\ait{\aitvar+1}$ in our case, for executing the first matrix vector multiplication. All intermediate results toward the optimum must be available. The implementation of the backward mode \eqref{eq:alg-iter-reverse} is shown in Algorithm~\ref{alg:bilevel-abstract-alg-reverse}.
\begin{figure*}[t]
\centering
\fbox{
\begin{minipage}{0.7\textwidth}
\begin{algorithm}\label{alg:bilevel-abstract-alg-reverse}
\underline{Derivative of an abstract algorithm}
\begin{itemize}
\item \emph{Assumptions}: $\Alg$ is totally differentiable.
\item \emph{Initialization at $\aitvar+1$}: 
    \[
      z^\ait{\aitvar+1}:= \Bigg(\frac{\partial \Loss}{\partial x} (x^*(\theta), \theta)\Bigg)^\top \in\R^\dimOpt \quad\text{and}\quad w^\ait{\aitvar+1}:= 0\in\R^\dimP
     \]
\item \emph{Iterations $(n\ge 0)$}: Update
    \[
      \begin{split}
        &\texttt{for $\aitvar$ to $0$}: \\
        &\left\lfloor
          \begin{split}
            w^\ait\aitvar =&\  w^\ait{\aitvar +1} + \Bigg(\frac{\partial \Alg}{\partial \theta} (x^\ait\aitvar,\theta)\Bigg)^\top z^\ait{\aitvar+1} \\
            z^\ait\aitvar =&\ \Bigg(\frac{\partial\Alg}{\partial x^\ait\aitvar}(x^\ait\aitvar,\theta)\Bigg)^\top z^\ait{\aitvar+1} 
          \end{split}
          \right.
      \end{split}
    \]
\item \emph{Final derivative of $\Loss$ in \eqref{eq:bilevel-abstract-alg} wrt. $\theta$}: 
    \[
      \frac{d\Loss}{d\theta} (x^*(\theta),\theta) = (w^\ait0)^\top + \frac{\partial\Loss}{\partial\theta}(x^*(\theta),\theta)\,.
    \]
  \end{itemize}
\end{algorithm}
\end{minipage}
}
\end{figure*}
This approach is quite expensive. But, for a reasonable number of iterations, it is still practical. It is still faster than the inversion of the Hessian matrix in Section~\ref{sec:derivative-impl-fun}; see \eqref{eq:full-derivative-impl-diff-param}. In the following section, we present approximations that reduce the cost significantly.

\subsection{Derivative of fixed point equations} \label{sec:derivative-fixed-point}

We generalize the result from Section~\ref{sec:derivative-impl-fun}, where the lower level problem of \eqref{eq:bilevel-general} is replaced by the first-order optimality condition of a smooth approximation. The idea is to consider a different optimality condition. A point is optimal, if it satisfies the fixed point equation of an algorithm $\map\Alg{X\times \R^\dimP}{X}$ solving the original lower level problem, i.e., we address the bilevel problem:
\begin{equation}\label{eq:bilevel-abstract-fp}
\begin{split}
  &\ \min_\theta\ \Loss(x^*(\theta),\theta) \\
  &\ \st\ x^*(\theta) = \Alg (x^*(\theta), \theta) \,,
\end{split}
\end{equation}
where $X\subset \R^\dimOpt$ is as in Section~\ref{sec:derivative-iterative-alg} and we have a fixed point $x^*$. This approach is more general than the one in Section~\ref{sec:derivative-impl-fun}, since we could actually first smoothly approximate the lower level problem and then consider the fixed point equation. For many algorithms  both approaches are equivalent, because optimization algorithms are often derived from the first-order optimality condition. 

Following the idea of Section~\ref{sec:derivative-iterative-alg}, we can consider a differentiable fixed point equation without the lower level problem to be differentiable. An algorithm that has a differentiable update rule yields a differentiable fixed point equation. 

Assume that $(x^*,\theta)$ solves the fixed point equation.  By differentiating the fixed point equation, we obtain
\[
  \frac{d x}{d \theta}(\theta) = \frac{\partial \Alg}{\partial x}(x^*(\theta),\theta) \frac{d x}{d \theta}(\theta)+ \frac{\partial \Alg}{\partial \theta}(x^*(\theta),\theta) \,,
\]
which can be rearranged to yield
\begin{equation} \label{eq:alg-impl-diff}
  \frac{d x}{d \theta}(\theta) = \Bigg(I - \frac{\partial \Alg}{\partial x}(x^*(\theta),\theta)\Bigg)^{-1} \frac{\partial \Alg}{\partial \theta}(x^*(\theta),\theta) \,.
\end{equation}
Assuming the spectral radius of $({\partial \Alg}/{\partial x})(x^*(\theta),\theta)$ is smaller than $1$, we can approximate the inversion using the geometric series:
\[
  \frac{d x}{d \theta}(\theta) = \sum_{\aitvar=0}^{\infty} \Bigg(\frac{\partial \Alg}{\partial x}(x^*(\theta),\theta)\Bigg)^{\aitvar} \frac{\partial \Alg}{\partial \theta}(x^*(\theta),\theta) \,,
\]
where $(({\partial \Alg}/{\partial x})(x^*(\theta),\theta))^\aitvar$ means the $\aitvar$-fold matrix product with itself. Let us approximate this term with a finite summation of $0,\ldots, \aitvar_0$. Then by a simple rearrangement, for $\xi\in\R^\dimOpt$,  we have (by abbreviating $({\partial \Alg}/{\partial x})(x^*(\theta),\theta)$ by ${\partial \Alg}/{\partial x}$; the same for ${\partial \Alg}/{\partial \theta}$):
\[
  \begin{split}
   \xi^\top \frac{d x}{d \theta}(\theta)\approx &\  
    \xi^\top \sum_{\aitvar=0}^{\aitvar_0} \Bigg(\frac{\partial \Alg}{\partial x}\Bigg)^{\aitvar} \frac{\partial \Alg}{\partial \theta} \\
    =&\ \xi^\top\frac{\partial \Alg}{\partial x}\left( \frac{\partial \Alg}{\partial \theta} + \frac{\partial \Alg}{\partial x} \left( \frac{\partial \Alg}{\partial \theta} + \ldots \right)\right) \\
    =&\ \xi^\top \left[ \frac{\partial\Alg}{\partial x^\ait{\aitvar_0}} \frac{dx^\ait{\aitvar_0}}{d\theta} \right]+ \xi^\top \frac{\partial \Alg}{\partial \theta} \,.
  \end{split}
\]
The difference between the last line in this equation and \eqref{eq:alg-iter-fwd}~and~\eqref{eq:alg-iter-reverse} is the evaluation point of the terms. While in \eqref{eq:alg-iter-fwd}~and~\eqref{eq:alg-iter-reverse} the terms for ${dx^\ait{\aitvar+1}}/{d\theta}$ are evaluated at $(x^\ait\aitvar(\theta),\theta)$, here, all terms are evaluated at $(x^*(\theta),\theta)$. Although the condition on the spectral radius is rarely met in practice, this approximation works well empirically and needs to store only the optimum of the algorithm. This leads to an immense reduction of the memory requirements.

\subsection{Weak differentiation of iterative algorithms} \label{sec:gen-derivative-alg}

The approach in \cite{DVFP14} also considers an algorithm replacing the non-smooth lower level problem. Their underlying methodology, however, is based on weak differentiability, which can be guaranteed for Lipschitz continuous mappings thanks to Rademacher's theorem. If all iteration mappings are Lipschitz continuous with respect to the iteration variable and the parameter $\theta$, weak differentiability follows from the chain rule for Lipschitz mappings \cite[Theorem 4]{EG92}. For details, we refer to \cite{DVFP14}, in particular Section 4.

\section{Explicit derivatives for exemplary algorithms} \label{sec:reverse-specific-algs}

The framework of Bregman proximity functions is key for the idea to approximate a non-smooth optimization problem by an algorithm with smooth update mappings. In this section, we instantiate two such algorithms. Details and examples of Bregman proximity functions are postponed to Section~\ref{sec:intro-bregman}. For understanding this section, it suffices to know that $D_\psi(x,\bar x)$ provides a distance measure between two points $x$ and $\bar x$, and it can be used to define a Bregman proximity operator $\prox[\psi]$ which generalizes the common proximity operator that is based on the Euclidean distance.

\subsection{Derivative of forward--backward splitting} \label{sec:FB-reverse}

Let us consider  forward--backward splitting \cite{LM79,Passty79} with Bregman proximity function $D_\psi$ (e.g. \cite{BT03}). It applies to minimization problems of the form
\[
  \min_{x\in\R^\dimOpt}\, f(x) + g(x) \,,
\]
where $\map{f}{\R^\dimOpt}{\R}$ is a continuously differentiable, convex function with Lipschitz continuous gradient and $\map{g}{\R^\dimOpt}{\eR}$ is a proper, lower semi-continuous, convex function with a (Bregman) proximity operator that is easy to evaluate. The update rule of the forward--backward splitting we consider is:
\begin{equation} \label{eq:FBS-update-abstract}
  \begin{split}
  x^\ait{\aitvar+1} =&\  \arg\min_{x\in\R^\dimOpt}\, g(x;\theta) + f(x^\ait\aitvar;\theta) \\
                     &\ + \scal{\nabla f(x^\ait\aitvar;\theta)}{x-x^\ait\aitvar} + \frac 1\alpha D_\psi(x, x^\ait\aitvar) \\
  =:&\ \prox[\psi]_{\alpha g}\Big(\nabla\psi(x^\ait\aitvar) - \alpha \nabla f(x^\ait\aitvar;\theta) ;\theta\Big) \\
  =:&\  \prox[\psi]_{\alpha g}\Big(y^\ait\aitvar(x^\ait\aitvar;\theta) ;\theta\Big) \,,
  \end{split}
\end{equation}
where we denote $y^\ait\aitvar(x^\ait\aitvar;\theta) := \nabla\psi(x^\ait\aitvar) - \alpha \nabla f(x^\ait\aitvar;\theta)$, the intermediate result after the forward step. The implementation of the reverse mode for determining the derivative of the solution map of the lower level problem with respect to $\theta$ is given in Algorithm~\ref{alg:bilevel-fb-alg-reverse}.
\begin{figure*}[t]
\centering
\fbox{
\begin{minipage}{0.7\textwidth}
\begin{algorithm}\label{alg:bilevel-fb-alg-reverse}
\underline{Derivative of a forward--backward splitting algorithm}
\begin{itemize}
\item \emph{Assumptions}: $\prox[\psi]_{\alpha g}$ and $\mathrm{id} + \alpha \nabla f$ are totally differentiable.
\item \emph{Initialization at $\aitvar+1$}: 
    \[
      z^\ait{\aitvar+1}:= \Bigg(\frac{\partial \Loss}{\partial x} (x^*(\theta), \theta)\Bigg)^\top \in\R^\dimOpt \quad\text{and}\quad w^\ait{\aitvar+1}:= 0\in\R^\dimP
     \]
\item \emph{Iterations $(n\ge 0)$}: Update (where derivatives of $\prox[\psi]_{\alpha g}$ are evaluated at $(y^\ait\aitvar,\theta)$ and derivatives of $\nabla f$ at $(x^\ait\aitvar;\theta)$)
    \[
      \begin{split}
        &\texttt{for $\aitvar$ to $0$}: \\
        &\left\lfloor
          \begin{split}
            w^\ait\aitvar =&\  w^\ait{\aitvar+1} + \left(\Bigg(\frac{\partial \prox[\psi]_{\alpha g}}{\partial \theta} \Bigg)^\top + \Bigg(- \alpha \frac{\partial (\nabla f)}{\partial \theta}\Bigg)^\top \Bigg(\frac{\partial \prox[\psi]_{\alpha g}}{\partial y}\Bigg)^\top \right) z^\ait{\aitvar+1} \\
            z^\ait\aitvar =&\ \Bigg(\mathrm{id} - \alpha \frac{\partial (\nabla f)}{\partial x}\Bigg)^\top \Bigg(\frac{\partial \prox[\psi]_{\alpha g}}{\partial y}\Bigg)^\top z^\ait{\aitvar+1} 
          \end{split}
          \right.
      \end{split}
    \]
\item \emph{Final derivative of $\Loss$ in \eqref{eq:bilevel-abstract-alg} wrt. $\theta$}: 
    \[
      \frac{d\Loss}{d\theta} (x^*(\theta),\theta) = (w^\ait0)^\top + \frac{\partial\Loss}{\partial\theta}(x^*(\theta),\theta)\,.
    \]
  \end{itemize}
\end{algorithm}
\end{minipage}
}
\end{figure*}

\subsection{Derivative of primal--dual splitting} \label{sec:PD-reverse}

Since the primal--dual algorithm with Bregman proximity functions from \cite{CP15} provides us with a flexible tool, we specify the implementation of the reverse mode for this algorithm. It applies to the convex--concave saddle-point problem
\[
  \min_x\max_y \scal{Kx}{y} + f(x) + g(x) - h^*(y)\,,
\]
which is derived from $\min_x f(x) + g(x) + h(Kx)$, where $f$ is convex and has a Lipschitz continuous gradient and $g,h$ are proper, lower semi-continuous convex functions with simple proximity operator for $g$ and for the convex conjugate $h^*$.

Let the forward iteration of the primal--dual algorithm with variables $x^\ait\aitvar=(u^\ait\aitvar,p^\ait\aitvar)\in\R^{\dimOpt_u+\dimOpt_p}$ be given as
\begin{equation}\label{eq:PD-abstract}
  \begin{split}
    u^\ait{\aitvar+1} =&\ \PDx(u^\ait\aitvar,p^\ait\aitvar,\theta) \\
                      :=&\ \arg\min_u\, \scal{\nabla f(u^\ait\aitvar)}{u-u^\ait\aitvar} + g(u) \\
                      &\ \qquad\quad + \scal{Ku}{p^\ait\aitvar} + \tfrac 1\tau D_u(u,u^\ait\aitvar) \\
    p^\ait{\aitvar+1} =&\ \PDy(2u^\ait{\aitvar+1}-u^\ait\aitvar, p^\ait\aitvar, \theta) \\
                      :=&\ \arg\min_p\, h^*(p) - \scal{K(2u^\ait{\aitvar+1}-u^\ait\aitvar)}{p} \\
                      &\ \qquad\quad + \tfrac 1\sigma D_p(p,p^\ait\aitvar) \,,
  \end{split}
\end{equation}
where $f,g,h,K$ can depend on $\theta$. The step size parameter $\tau$ and $\sigma$ must be chosen according to $(\tau^{-1} - L_f) \sigma^{-1} \geq L^2$ where $L=\norm{K}$ is the operator norm of $K$ and $L_f$ is the Lipschitz constant of $\nabla f$. 

\begin{figure*}[t]
  \begin{center}
  \newcommand{\thisTikzScaling}{1.0}
\newcommand{\DERIV}[2]{\partial_{#2}#1}
\newcommand{\TDERIV}[2]{\frac{d#1}{d#2}}
\begin{tikzpicture}[scale=\thisTikzScaling]

\begin{scope}
\matrix (m) [matrix of math nodes,row sep=3.5em,column sep=8em,minimum width=2em]
{
        & \TDERIV{u^*}{\theta}     &                                      &         \\
d\theta & \TDERIV{u^\ait\aitvar}{\theta}       & \TDERIV{p^\ait\aitvar}{\theta}       &         \\
d\theta & \TDERIV{u^{\ait{\aitvar-1}}}{\theta} & \TDERIV{p^{\ait{\aitvar-1}}}{\theta} & d\theta \\
d\theta & \TDERIV{u^{\ait{\aitvar-2}}}{\theta} & \TDERIV{p^{\ait{\aitvar-2}}}{\theta} & d\theta \\};

\tiny

\path[-stealth]
    (m-1-2) edge[gray] node [black,left] {$\DERIV{\PDx^\ait{\aitvar}}{u}$} (m-2-2)
            edge[gray,bend left=15] node [black,near end] {$\DERIV{\PDx^\ait{\aitvar}}{p}$} (m-2-3)
            edge[gray,bend right=15] node[black,left] {$\DERIV{\PDx^\ait{\aitvar}}{\theta}$} (m-2-1)
    (m-2-2) edge[gray] node [black,left] {$\DERIV{\PDx^{\ait{\aitvar-1}}}{u}$} (m-3-2)
            edge[gray,bend left=15] node [black,near end] {$\DERIV{\PDx^{\ait{\aitvar-1}}}{p}$} (m-3-3)
            edge[gray,bend right=15] node[black,left] {$\DERIV{\PDx^{\ait{\aitvar-1}}}{\theta}$} (m-3-1)
    (m-3-2) edge[gray] node [black,left] {$\DERIV{\PDx^{\ait{\aitvar-2}}}{u}$} (m-4-2)
            edge[gray,bend left=15] node [black,near end] {$\DERIV{\PDx^{\ait{\aitvar-2}}}{p}$} (m-4-3)
            edge[gray,bend right=15] node[black,left] {$\DERIV{\PDx^{\ait{\aitvar-2}}}{\theta}$} (m-4-1);

\path[-stealth]
    (m-2-3) edge[gray] node [black,right] {$\DERIV{\PDy^{\ait{\aitvar-1}}}{p}$} (m-3-3)
            edge[gray,bend left=15] node [black,right] {$\DERIV{\PDy^{\ait{\aitvar-1}}}{\theta}$} (m-3-4)
            edge[gray] node [black,above] {$2\DERIV{\PDy^{\ait{\aitvar-1}}}{u}$} (m-2-2)
            edge[gray,bend right=15] node [black,near end] {$-\DERIV{\PDy^{\ait{\aitvar-1}}}{u}$} (m-3-2)
    (m-3-3) edge[gray] node [black,right] {$\DERIV{\PDy^{\ait{\aitvar-2}}}{p}$} (m-4-3)
            edge[gray,bend left=15] node [black,right] {$\DERIV{\PDy^{\ait{\aitvar-2}}}{\theta}$} (m-4-4)
            edge[gray] node [black,above] {$2\DERIV{\PDy^{\ait{\aitvar-2}}}{u}$} (m-3-2)
            edge[gray,bend right=15] node [black,near end] {$-\DERIV{\PDy^{\ait{\aitvar-2}}}{u}$} (m-4-2)
    (m-4-3) edge[gray] node [black,above] {$2\DERIV{\PDy^{\ait{\aitvar-3}}}{u}$} (m-4-2);

\foreach \x in {1,2,3,4}
  \node[below=0.5cm] at (m-4-\x) {$\vdots$};

\end{scope}
\end{tikzpicture}
  \end{center}
  \caption{\label{fig:back-prop-pd}The graph shows how the information is backprogated to estimate the derivatives in Algorithm~\ref{alg:bilevel-pd-alg-reverse}. The derivatives at the nodes show what derivative is to be evaluated from this point downwards through the graph. The edges represent multiplicative (transposed) factors. The final derivative is the sum over all leaf nodes.}
\end{figure*}
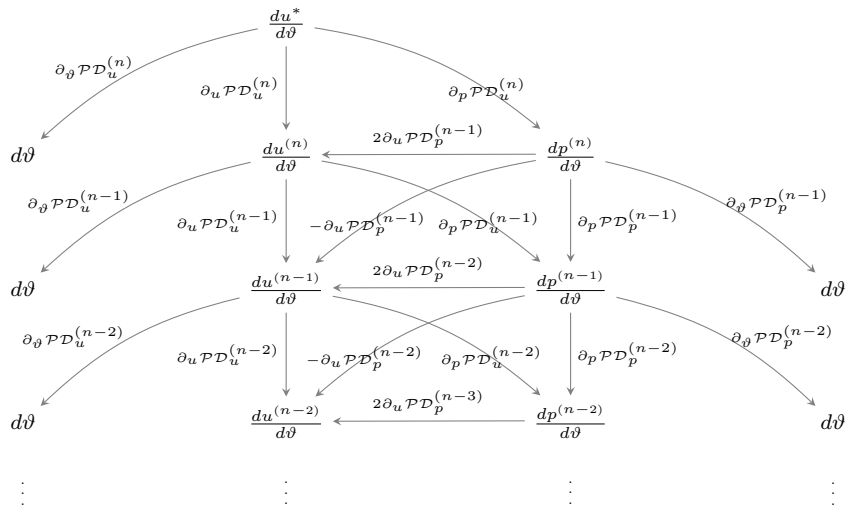

To illustrate the application of the chain rule throughout the primal--dual algorithm, we show a graphical representation of the information flow in Figure~\ref{fig:back-prop-pd}, where we use the following abbreviations (analogously for $\PDy$): 
\begin{gather*}
  \PDx^\ait\aitvar := \PDx(u^\ait\aitvar,p^\ait\aitvar,\theta)\,; \\
  \PDy^\ait\aitvar := \PDy(2u^\ait{\aitvar+1}-u^\ait{\aitvar},p^\ait\aitvar,\theta)\,; \\
  \partial_u \PDx := \frac{\partial \PDx}{\partial u}\,;
  \partial_p \PDx := \frac{\partial \PDx}{\partial p}\,;
  \partial_\theta \PDx := \frac{\partial \PDx}{\partial \theta} \,.
\end{gather*}
\begin{remark}
  In Section~\ref{sec:FB-reverse}, we evaluated the forward and the backward step separately using the chain rule. Of course, this could be done here as well.
\end{remark}
Based on this graphical representation, it is easy to derive Algorithm~\ref{alg:bilevel-pd-alg-reverse}.
\begin{figure*}[t]
\centering
\fbox{
\begin{minipage}{0.7\textwidth}
\begin{algorithm}\label{alg:bilevel-pd-alg-reverse}
\underline{Derivative of a primal--dual algorithm}
\begin{itemize}
\item \emph{Assumptions}: $\PDx$ and $\PDy$ are totally differentiable.
\item \emph{Initialization at $\aitvar+1$}: 
    \begin{gather*}
      z^\ait{\aitvar+1}:= \Bigg(\frac{\partial \Loss}{\partial u} (u^*(\theta), \theta)\Bigg)^\top \in\R^{\dimOpt_u}\,, \quad q^\ait{\aitvar+1} := 0 \in \R^{\dimOpt_p} \\
      \text{and}\quad w^\ait{\aitvar+1}:= 0\in\R^\dimP
     \end{gather*}
\item \emph{Iterations $(n\ge 0)$}: Update
    \[
      \begin{split}
        &\texttt{for $\aitvar$ to $0$}: \\
        &\left\lfloor
          \begin{split}
            w^\ait\aitvar =&\  w^\ait{\aitvar +1} + \Bigg(\frac{\partial \PDx^\ait\aitvar}{\partial \theta}\Bigg)^\top z^\ait{\aitvar+1} + \Bigg(\frac{\partial \PDy^\ait\aitvar}{\partial \theta}\Bigg)^\top q^\ait{\aitvar+1}\\
            q^\ait\aitvar =&\ \Bigg(\frac{\partial\PDx^\ait\aitvar}{\partial p}\Bigg)^\top z^\ait{\aitvar+1} + \Bigg(\frac{\partial\PDy^\ait\aitvar}{\partial p}\Bigg)^\top q^\ait{\aitvar+1} \\
            z^\ait\aitvar =&\ \Bigg(\frac{\partial\PDx^\ait\aitvar}{\partial u}\Bigg)^\top z^\ait{\aitvar+1} + 2 \Bigg(\frac{\partial\PDy^\ait{\aitvar-1}}{\partial u}\Bigg)^\top q^\ait{\aitvar} - \Bigg(\frac{\partial\PDy^\ait\aitvar}{\partial u}\Bigg)^\top q^\ait{\aitvar+1} \\
          \end{split}
          \right.
      \end{split}
    \]
\item \emph{Final derivative of $\Loss$ in \eqref{eq:bilevel-abstract-alg} with $\Alg=(\PDx,\PDy)$ wrt. $\theta$}: 
    \[
      \frac{d\Loss}{d\theta} (u^*(\theta),\theta) = (w^\ait0)^\top + \frac{\partial\Loss}{\partial\theta}(u^*(\theta),\theta)\,.
    \]
  \end{itemize}
\end{algorithm}
\end{minipage}
}
\end{figure*}

A running average is used to implement the ergodic primal--dual algorithm whose output is the average of all iterates, i.e., $u^* = \frac 1{\aitvar +1}\sum_{i=0}^{\aitvar} u^\ait i$: denote $s_u^\ait{\aitvar}:= \frac 1{\aitvar+1}\sum_{i=0}^{\aitvar} u^\ait i$, then $s_u^{\ait{\aitvar+1}} = \frac{1}{\aitvar+2} u^\ait{\aitvar+1} + \frac{\aitvar+1}{\aitvar+2}s_u^\ait{\aitvar}$. Since the derivative is a linear operator, we can estimate the derivative for the ergodic primal--dual sequence by averaging all $w^\ait\aitvar$. These can be computed as a running average in the loop of Algorithm~\ref{alg:bilevel-pd-alg-reverse}.

\section{\enquote{Smoothing} using Bregman proximity} \label{sec:ex-technique}

Splitting based techniques like those in Section~\ref{sec:reverse-specific-algs} usually handle non-smooth terms in the objective function via a (non-linear/Bregman) proximal step. Convex conjugation makes terms in the objective amenable for simple and differentiable proximal mappings. Adding the possibility of considering a primal, primal--dual, or dual formulation yields many examples of practical interest.

In the following, we introduce the class of Bregman functions that can be used in combination with the algorithms in Section~\ref{sec:reverse-specific-algs}. Then, we discuss a few examples that allow the reformulation of several non-smooth terms arising in applications.

\subsection{Bregman proximity functions} \label{sec:intro-bregman}

We consider Bregman proximity functions \cite{Bregman67} with the following properties: Let $\map\psi {\R^\dimOpt}\eR$ be a 1-convex function with respect to the Euclidean norm, i.e., it is strongly convex with modulus $1$, and denote its domain by $X:=\dom \psi$. We assume that $\psi$ is continuously differentiable on the interior of its domain $\int X$ and continuous on its closure $\cl X$. 

Then, $\psi$ generates a Bregman proximity function $\map{D_\psi}{X \times \int X}{\R}$ by
\begin{equation} \label{eq:def-bregman-gen}
  D_\psi(x,\bar x) := \psi(x) - \psi(\bar x) - \scal{\nabla\psi(\bar x)}{x-\bar x}\,.
\end{equation}
For a sequence $\seq[n\in\N]{x^n}$ converging to $x\in X$, we require that $\lim_{n\to\infty} D_\psi(x,x^n) = 0$. The 1-convexity of $\psi$ implies that the Bregman function satisfies the inequality
\[
  D_\psi(x,\bar x) \geq \frac 12 \norm{x-\bar x}^2\,, \quad \forall x\in X,\; \bar x \in \int X\,.
\]
These are the kind of Bregman proximity functions considered in \cite{CP15}. Obviously $\psi(x)=\frac12 \norm{x}^2$ corresponds to $D_\psi(x,\bar x) = \frac 12 \norm{x-\bar x}^2$. 

In iterative algorithms, the Bregman proximity function is used via the proximity operator for a proper, lower semi-continuous, convex function $\map{g}{X}{\eR}$
\begin{equation} \label{eq:Bregman-prox-min}
  \prox[\psi]_{\alpha g} (\bar x) := \arg\min_{x\in X}\, \alpha g(x) + D_\psi(x,\bar x) \,,
\end{equation}
where we define $\prox_{\alpha g}:= \prox[\frac 12\norm{\cdot}^2]_{\alpha g}$.\\

There are two kinds of Bregman proximity functions: \titemi The function $\nabla \psi$ can be continuously extended to $X$, i.e., $D_\psi$ can be defined on $X\times X$, and \titemii $\psi$ is differentiable on $\int X$ (i.e. $\nabla \psi$ cannot necessarily be extended to $\cl X$). In this case $D_\psi(x,\bar x)$ makes sense only on $X\times\int X$ and we must assure that $\prox[\psi]_{\alpha g}(\bar x)\in\int X$ for any $\bar x\in \int X$. For this, we need to assume that $\norm{\nabla\psi(x)}\to\infty$ whenever $x$ approaches a boundary point $\bdry X := \cl X \smallsetminus\int X$ (which is sometimes referred to as $\psi$ being essentially smooth \cite{Rock70}). 

While solutions of the proximity operator for the first class can lie on the boundary $\bdry X$, this is not possible for the second class; boundary points can be reached only in the limit when the proximity operator is applied sequentially. Moreover, for $\bar x\in\bdry X$, \eqref{eq:def-bregman-gen} would imply that, unless $x=\bar x$, the Bregman distance is $+\infty$ for any $x$, which can be represented by $\delta_{[x=\bar x]}(x)$. This means $\bar x\in\bdry X$ is always a fixed point of this Bregman proximity operator. This precludes application of the fixed-point approach from Section~\ref{sec:derivative-fixed-point}.

\subsection{Examples of Bregman functions} \label{subsec:Bregman-examples}

Since Bregman proximity functions play a key role in this paper, we consider a few examples.
\begin{example} \label{ex:Bregman-Euclid}
    The Euclidean length $\psi(x) = \frac 12\norm[2]{x}^2$ is continuously differentiable on the whole space $\R^\dimOpt$, and therefore, belongs to class \titemi of Bregman proximity functions. 
\end{example}
\begin{example} \label{ex:Bregman-interval}
    The Bregman proximity function generated by $\psi(x) = \frac 12((x+1)\log(x+1) + (1-x)\log(1-x))$ is defined on the interval $(-1,1)$ and can be continuously extended to $[-1,1]$, and is continuously differentiable on $(-1,1)$ with $\abs{\psi^\prime(x)} \to \infty$ when $x\to \pm 1$. It is 1-strongly convex. 
\end{example}
\begin{example}\label{ex:bregman-entropy}
    The entropy function $\psi(x) = x\log(x)$, which can be continuously extended to $[x\geq 0]$, is continuously differentiable on $[x>0]$ with derivative $\psi^\prime(x) = \log(x)+1$.  The derivative cannot be continuously extended to $x=0$. For $x\to 0$ we have $\abs{\psi^\prime(x)}\to +\infty$. Unfortunately, this function is not even 1-strongly convex on $[x\geq 0]$. However, the function $a x\log(x)$ is 1-strongly convex when restricted to a bounded subset $[0,1/a]$, $a>0$. For $a=1$, the Bregman function $D_\psi(x,\bar x)=x(\log(x)-\log(\bar x))-(x-\bar x)$ is generated.
\end{example}
\begin{example}
  The entropy function can also be used in higher dimensions. Unfortunately, it is hard to assert a simple evaluation of an associated proximity mapping in this case. Consider a polyhedral set $\emptyset \neq X\in\R^\dimOpt$ given by 
  \[
    \begin{split}
    X =&\  \set{x\in\R^N\vert\, \forall i=1,\ldots,M\colon \scal{a_i}{x} \leq b_i}  \\
    =&\  \bigcap_{i=1}^M \set{x \in\R^\dimOpt\vert\, \scal{a_i}{x} \leq b_i}
    \end{split}
  \]
  for vectors $0\neq a_i\in \R^{\dimOpt}$, and $b_i\in\R^M$, $i=1,\ldots,M$. Then, the generating function 
  \[
    \psi (x) = \sum_{i=1}^M (b_i-\scal{a_i}{x})\log(b_i-\scal{a_i}x) 
  \]
  is designed such that for any point $\bar x \in \int X$ any other point $x\not\in X$ is \enquote{moved infinitly far away} with respect to the Bregman distance $D_\psi(x,\bar x)$. Therefore $\norm{\nabla \psi(x)} \to \infty$ for $x$ tends towards a point on the boundary $\bdry X$. Nevertheless, $\psi$ is continuous on $X$ and strongly convex, if $X$ is bounded.
\end{example}

\subsection{Examples of smooth Bregman proximity operators}

The Bregman proximity functions that we presented are particularly interesting if the evaluation of the proximal mapping \eqref{eq:Bregman-prox-min} is a constrained minimization problem, i.e. the involved function $g$ in $\prox[\psi]_{g}$ is extended-valued and $+\infty$ outside the constraint (closed) convex set $X\subset\R^\dimOpt$. The Bregman function can replace or simplify the constraint set. In the following, we consider a few examples of practical interest. The class of functions that are amenable to our approach can be broadened significantly thanks to the concept of (convex) conjugation.

We consider a basic class of functions $g(x) = \scal{x}{c} + \ind{X}(x)$ for some $c\in\R^\dimOpt$. The associated (non-linear) proximity operator from \eqref{eq:Bregman-prox-min} is given by
\[
 \prox[\psi]_{\alpha g} (\bar x) = \arg\min_{x\in X}\, \alpha \scal{x}{c} + D_\psi(x,\bar x)  \,.
\]
The corresponding (necessary and sufficient) optimality condition, which has a unique solution, is
\[
  \begin{split}
 &\ 0 \in c + \nabla\psi(x) - \nabla\psi(\bar x) + \partial \ind{X}(x) \\
 \Leftrightarrow&\ \nabla\psi(\bar x) - c \in \nabla\psi(x) + \ncone X(x) \,,
 \end{split}
\]
where $\ncone X(x)$ denotes the normal cone at $x$ of the set $X$. Suppose $\bar x \in \int X$. If $\psi$ is chosen such that $\norm{\nabla \psi(x)} \to +\infty$ for $x \to \tilde x\in\bdry X$, then the solution of the proximal mapping is in $\int X$. Since $\ncone X(x)=0$ for $x\in\int X$, the optimality condition simplifies to 
\begin{equation} \label{eq:Bregman-prox-lin-constr}
  \nabla\psi(\bar x) - c = \nabla\psi(x) \,,
\end{equation}
i.e. the constraint is implicitly taken care of by the Bregman proximity function. Summarizing, the goal of our approach  (for this basic function $g$) consists of determining $\psi$, respectively $D_\psi$, such that
\begin{itemize}
  \item the constraint set can be handled implicitly, 
  \item \eqref{eq:Bregman-prox-lin-constr} can be solved efficiently (possibly in closed form),
  \item and the solution function of \eqref{eq:Bregman-prox-lin-constr}, which yields the solution of \eqref{eq:Bregman-prox-lin-constr} for a given $\bar x$, is required to be differentiable wrt. $x$ and $\theta$, where possibly $c = c(\theta)$.
\end{itemize}
\begin{example}
    For a linear function $g(x) = \scal cx + \ind{[x\geq 0]}(x)$ the entropy function from Example~\ref{ex:bregman-entropy} can be summed-up for each coordinate to remove the non-negativity constraint. The proximity operator reads: 
    \[
      \Bigg(\prox[\sum_j x_j\log x_j]_{\alpha g}(\bar x)\Bigg)_i = \bar x_i \exp(-\alpha c_i) \,.
    \]
    A closer look at the iterations of the forward--backward splitting (FBS) algorithm \eqref{eq:FBS-update-abstract} reveals that such a function $g$ arises  with $c=\nabla f(\bar x)$, i.e. in the iterations of FBS for the minimization of
    \[
      \min_{x\in \R^\dimOpt}\, f(x)  + \ind{[x\geq 0]}(x) \,.
    \]
    A particular instance of this problem is the \emph{non-negative least squares problem}, i.e. $f(x) = \frac 12 \norm[2]{Ax - b}^2$ with a matrix $A$ and a vector $b$.
\end{example}

\begin{example} \label{ex:Bregman-prox-lin-simplex}
 The most frequent application of the entropy-prox is for the minimization of a linear function $g(x) = \scal{c}{x}$ over the unit simplex in $\R^\dimOpt$. Since the entropy function restricts the solution of the proximity operator to the positive orthant, projecting a point $\bar x\in\R_+^\dimOpt$ onto the unit simplex $\set{x\in\R^\dimOpt\vert\, \sum_{i=1}^\dimOpt x_i = 1\text{ and } x_i \geq 0}$ reduces to the projection onto the affine subspace $\set{x\in\R^\dimOpt\vert\, \sum_{i=1}^\dimOpt x_i = 1}$, which can be given in closed-form, i.e.,
    \[
      \Bigg( \prox[\sum_j x_j\log x_j]_{\alpha g} (\bar x)\Bigg)_i = \frac{\bar x_i\exp(-\alpha c_i)}{\sum_{j=1}^\dimOpt \bar x_j\exp(-\alpha c_j)} \,.
    \]
    This proximal problem arises for example in the \emph{multi-label segmentation problem} in Section~\ref{subsec:potts-model} or in \emph{Matrix games} (see \cite[Section 7.1]{CP15}).
\end{example}

\begin{example} \label{ex:Bregman-prox-lin-interval}
  For the function $g(x) = \scal cx + \ind{[-1\leq x\leq 1]}(x)$ the Bregman function from Example~\ref{ex:Bregman-interval} reduces the minimization problem in the proximal mapping to an unconstrained problem. The proximal mapping with $\psi(x) = \sum_i \frac 12((x_i+1)\log(x_i+1) + (1-x_i)\log(1-x_i))$ reads: 
    \[
      \Bigg(\prox[\psi]_{\alpha g}(\bar x)\Bigg)_i = \frac{\exp(-2\alpha c_i) - \frac{1-\bar x_i}{1+\bar x_i}}{\exp(-2\alpha c_i) + \frac{1-\bar x_i}{1+\bar x_i}} \,.
    \]
    Obviously, this example can be adjusted to any Cartesian product of interval constraints. The importance of this exemplary function $g$ becomes clear in the following.
\end{example}

Functions that are linear on a constraint set also arise when conjugate functions are considered. For instance the $\ell_1$-norm can be represented as 
\[
  \norm[1]{x} = \max_{y}\, \scal xy + \ind{[-1\leq y \leq 1]}(y)\,.
\]
In combination with the primal--dual (PD) algorithm \eqref{eq:PD-abstract}, this representation results in subproblems  of the type discussed in the preceding examples. From this perspective, optimization problems involving a linear operator $\opKnz$ and the $\ell_1$-norm $\norm[1]{\opKnz x}$ are also easy to address. 

This idea of conjugation can be put into a slightly larger framework, as the convex conjugate of any positively one-homogeneous proper, lsc, convex function is an indicator function of a closed convex set. Unfortunately, it is required that projecting onto such a set is easy (\enquote{prox-friendliness}). Therefore, the following example is restricted to the (additively) separable case.
\begin{example} \label{ex:Bregman-separable-conjugation}
  Let $g$ be an (additively) separable, positively one-homogeneous, proper, lsc, convex functions $g(x) = \sum_{i=1}^\dimOpt g_i(x_i)$. Thanks to its properties $g$ coincides with its bi-conjugate function $g^{**}$ and we can consider
  \[
      g(x) = g^{**}(x) = \sum_{i=1}^\dimOpt \max_{y_i}\, x_i y_i - \ind{Y_i}(y_i) \,,
  \]
  where $Y_i = [a_i,b_i]$ is a closed interval in $\R$. Again the dual update step of \eqref{eq:PD-abstract} involves problems such as in Example~\ref{ex:Bregman-prox-lin-interval} with $h^*(y) = \sum_i \ind{Y_i}(y_i)$.
\end{example}

\section{Toy example} \label{sec:toy-example}

The bilevel problem that we consider here is a parameter learning problem of a one dimensional non-negative least-squares problem:
\begin{equation} \label{eq:bilevel-nonneg-LS}
  \begin{split}
    \min_{\theta\in\R}&\; \frac12 (x^*(\theta) - \mathfrak g)^2 \\
      &\; \st\  x^*(\theta) = \arg\min_{x\in \R}\; \frac \lambda 2 (\theta x - b)^2 + \frac 12 x^2 + \delta_{[x\geq 0]}(x)\,,
  \end{split}
\end{equation}
where $\theta$ is the optimization variable of the bilevel problem, $b\in\R$ is the input of the least squares problem, and $\lambda$ is a positive weighting parameter. Given $\theta$ and $b$ the lower level problem solves the non-negative least squares problem. The squared Euclidean loss function in the upper level problem compares the output of the lower level problem for some $\theta$ and $b$ to the ground truth $\mathfrak g:=x^*(\theta^*)$, which is generated by solving the lower level problem with some predefined value $\theta^*$. The goal of the bilevel optimization problem is to find $\theta^*$ given $b$ and $\mathfrak g$. 
\begin{figure}[t]
\centering
\includegraphics[width=0.48\linewidth]{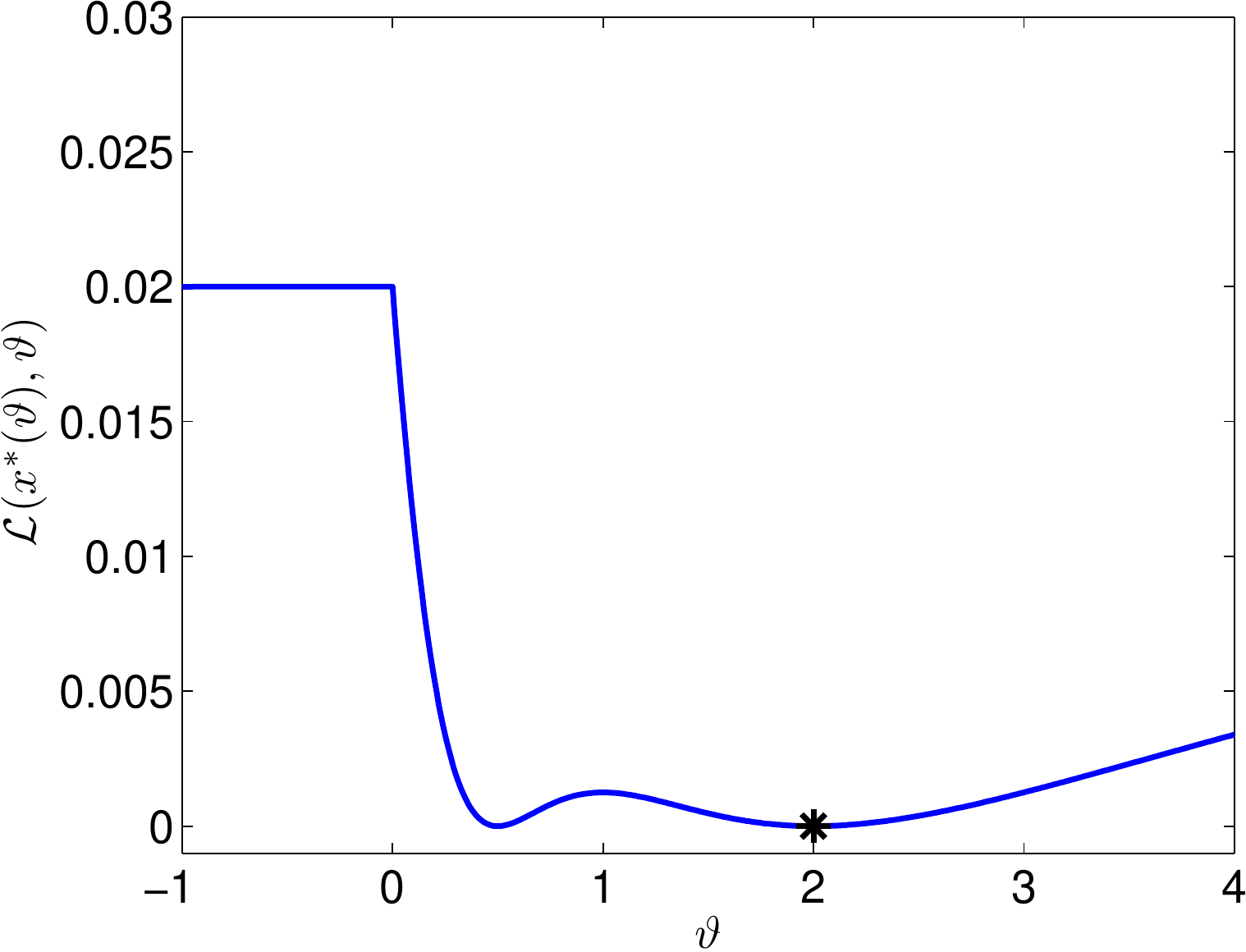}
\includegraphics[width=0.48\linewidth]{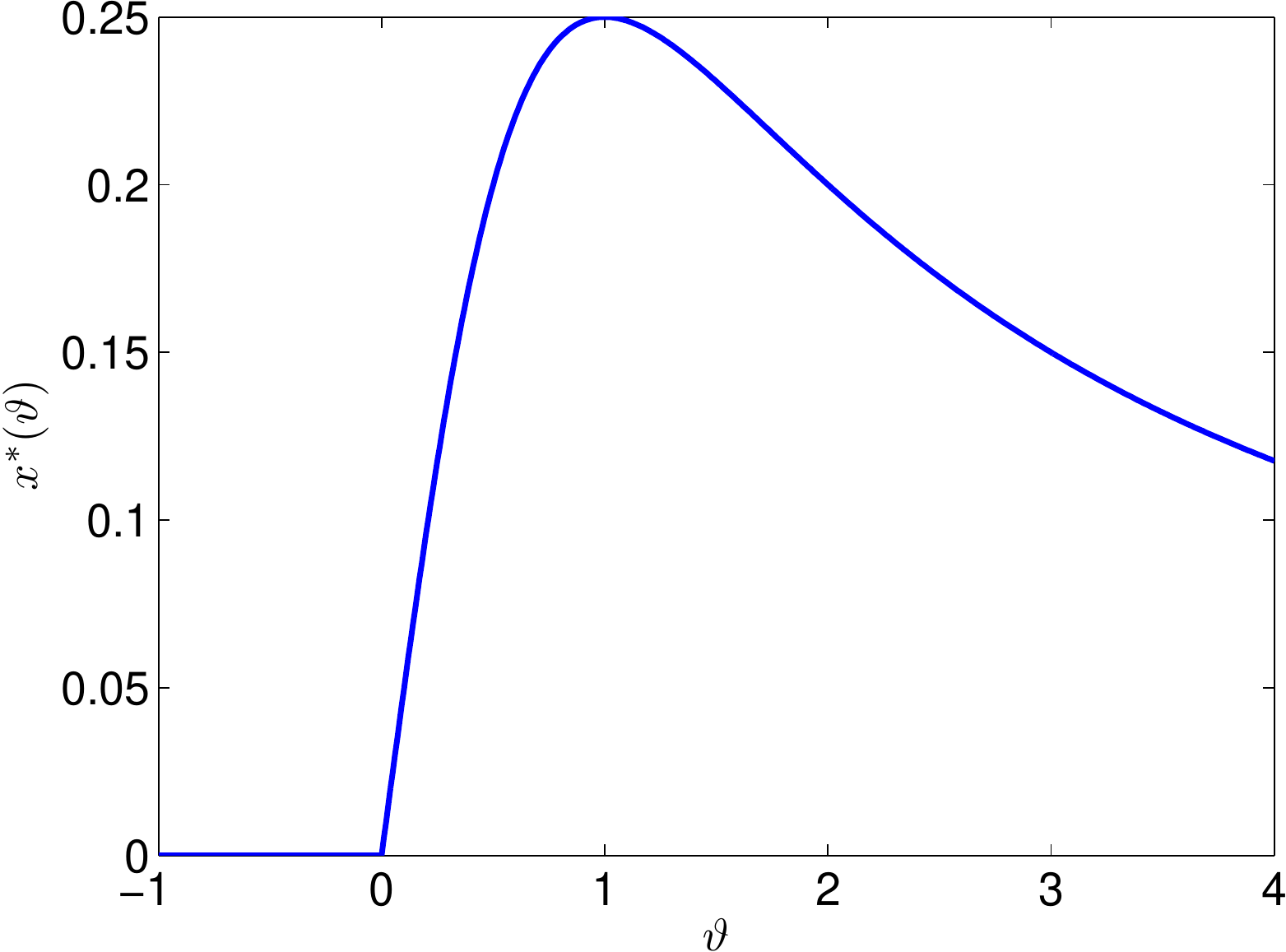}
\caption{\label{fig:example-nn-ls}Visualization of the loss function $\Loss(x(\theta),\theta)$ for the 1D example \eqref{eq:bilevel-nonneg-LS} on the left side. The optimum is marked with a black star. On the right hand side, the solution map of the lower level problem is shown.}
\end{figure}

The analytic solution of the lower level problem (the solution map) is 
\[
  x^*(\theta) = \max\Big(0, \frac{\lambda\theta b}{1+\lambda\theta^2}\Big)
\]
and is shown on the right hand side of Figure~\ref{fig:example-nn-ls}. It is obviously a non-smooth function with a non-differentiable point at $\theta=0$. Plugging the solution map into the upper level problem shows the actual objective to be minimized; see the left hand side of Figure~\ref{fig:example-nn-ls}. 

\subsection{Experimental setup} 

In the following experiments, we numerically explore the gradients computed with the proposed techniques. We do not consider the actual minimization of the bilevel problem. The computed gradients could be used by any first-order gradient based method.

\newcommand{\ImplDiff}{{\tt Smoothed-impl}\xspace}
\newcommand{\ProjGD}{\texttt{Proj.GD}\xspace}
\newcommand{\ProjGDb}{\texttt{Proj.GD2}\xspace}
\newcommand{\BregFB}{\texttt{Bregman-FB}\xspace}
\newcommand{\BregFBb}{\texttt{Bregman-FB2}\xspace}
\newcommand{\AlgImpl}{\texttt{Bregman-FB-impl}\xspace}

\paragraph{Analytic subdifferential.} For $\theta\neq 0$ the standard chain rule from calculus can be applied and we can directly write down the derivative of the whole problem, namely
\[
  \frac{d\Loss }{d \theta} (x(\theta)) = \frac{\lambda b(1-\lambda\theta^2)}{(1+\lambda\theta^2)^2} (x(\theta)-\mathfrak g)\,.
\]
For $\theta= 0$, we consider the derivative
\[
  \frac{d\Loss }{d \theta} (x(\theta)) = [0,\lambda b (x(0)-\mathfrak g)] \,,
\]
where $[0,\lambda b]$ is replaced by $[\lambda b, 0]$ if $\lambda b<0$.

\paragraph{Implicit differentiation approach Section~\ref{sec:derivative-impl-fun}.} In order to apply this technique, we must smooth the lower level problem. Since we want to avoid solutions $x^*(\theta)=0$, we introduce a log-barrier and replace the lower level problem by
\[
  f_\mu(x,\theta) := \frac \lambda2 (\theta x - b)^2 + \frac 12 x^2 - \mu \log(x)
\]
for some small $\mu>0$. Thus, we can drop the non-negativity constraint. To compute the gradient via the implicit differentiation formula \eqref{eq:full-derivative-impl-diff-param}, we minimize $f_\mu$ with respect to $x$ and compute the second derivatives (we abbreviate the $x$-derivative with $f_\mu^\prime$ and $\theta$-derivative with $\partial_\theta f_\mu$)
\begin{equation} \label{eq:example-nn-ls-derivatives}
\begin{split} 
  & f_\mu^\prime(x,\theta) = \lambda \theta(\theta x- b) + x - \frac \mu x\,; \\
  & f_\mu^{\prime\prime}(x,\theta) = \lambda \theta^2 + 1 + \frac \mu{x^2}\,; \\
  & \partial_\theta f^\prime_\mu (x,\theta) = 2\lambda \theta x -\lambda b \,.
\end{split}
\end{equation}
Then, \eqref{eq:full-derivative-impl-diff-param} yields
\[
    \frac{d\Loss }{d \theta} (x^*(\theta))= - (x(\theta)-\mathfrak g) (f_\mu^{\prime\prime}(x^*(\theta),\theta))^{-1} \partial_\theta f^\prime_\mu (x^*(\theta),\theta) \,.
\]
This approach is denoted \ImplDiff.

\paragraph{Algorithmic differentiation approach Section~\ref{sec:derivative-iterative-alg}.} We consider two algorithms: projected gradient descent and forward--backward splitting with Bregman proximity functions. Both algorithms are splitting methods that distribute the objective into a smooth function $f$ and a non-smooth function $g$, for our example it reads
\[
  f(x,\theta) = \frac \lambda 2 (\theta x - b)^2 + \frac 12 x^2\quadand g(x) = \delta_{[x\geq 0]}(x) \,.
\]
Projected gradient descent operates by a gradient descent step with respect to the smooth function $f$ followed by a projection onto the (convex) set $[x\geq 0]$:
\begin{equation} \label{eq:projGD-ex}
  \begin{split}
    x^\ait{\aitvar+1} =&\ \proj_{[x\geq 0]} (x^{\ait\aitvar} - \alpha f^\prime(x^\ait\aitvar,\theta) ) \\
                      =&\ \max(0, x^{\ait\aitvar} - \alpha f^\prime(x^\ait\aitvar) ) \,.
  \end{split}
\end{equation}
Note that the projection onto the convex set can also be interpreted as solving the proximity operator associated with the function $g$. 

The second algorithm is obtained by changing the distance function for evaluating the proximity operator to the Bregman distance from Example~\ref{ex:bregman-entropy}. It results in
\begin{equation} \label{eq:BregGD-ex}
  x^\ait{\aitvar+1} = x^n \exp(-\alpha f^\prime(x^\ait\aitvar,\theta)) \,.
\end{equation} 
As we assume that $x^0 \in [x> 0]$ the Bregman proximity function ensures that the solution stays in the feasible set. Thus, the back-projection can be dropped. \\

To apply Algorithm~\ref{alg:bilevel-abstract-alg-reverse} or~\ref{alg:bilevel-fb-alg-reverse}, we need the second derivatives of the update steps \eqref{eq:projGD-ex} and \eqref{eq:BregGD-ex}. The second derivatives of $f=f_\mu$ with $\mu=0$ are given in \eqref{eq:example-nn-ls-derivatives}. Although,  \eqref{eq:projGD-ex} is not differentiable, it is differentiable almost everywhere, and in the experiment, we formally applied the chain rule and assigned an arbitrary subgradient wherever it is not unique, i.e.,
\[
  \frac{\partial \proj_{[x\geq 0]}}{\partial x} (x,\theta) = \begin{cases}
    0 ,& \text{if } x< 0\,; \\
    1 ,& \text{if } x>0\,; \\
    [0,1] ,& \text{if } x=0\,;
  \end{cases}
\]
and $\frac{\partial \proj_{[x\geq 0]}}{\partial \theta} = 0$. This approach is denoted \ProjGD.

For \eqref{eq:BregGD-ex}, we use Algorithm~\ref{alg:bilevel-abstract-alg-reverse} and obtain\footnote{Note that we kept the order of the terms given by the chain rule, since for multi-dimensional problems the products are matrix products and are, in general, not commutative.}
\[
  \begin{split}
    \frac{\partial \Alg}{\partial \theta}(x^\ait\aitvar, \theta) =&\ -\alpha x \exp(-\alpha f^\prime(x^\ait\aitvar,\theta)) \frac{\partial f^{\prime}}{\partial \theta} (x^\ait\aitvar, \theta) \\
  \frac{\partial \Alg}{\partial x}(x^\ait\aitvar.\theta) =&\ \exp(-\alpha f^\prime(x^\ait\aitvar,\theta)) \\&\ -\alpha x^\ait\aitvar \exp(-\alpha f^\prime(x^\ait\aitvar,\theta))f^{\prime\prime} (x^\ait\aitvar, \theta) \,.
  \end{split} 
\]
This approach is denoted \BregFB.

\paragraph{Implicit differentiation of the fixed point equation approach from Section~\ref{sec:derivative-fixed-point}.} As explained above, direct differentiation of the fixed point equation of an algorithm implies two techniques. One is by applying Algorithm~\ref{alg:bilevel-abstract-alg-reverse} to \eqref{eq:BregGD-ex} but evaluating all derivatives at the optimum (denoted \BregFBb). The other is to do the numerical inversion as in \eqref{eq:alg-impl-diff} (denoted \AlgImpl).

\subsection{Analysis of the 1D example}

\begin{figure}[t]
\centering
\includegraphics[width=0.48\linewidth]{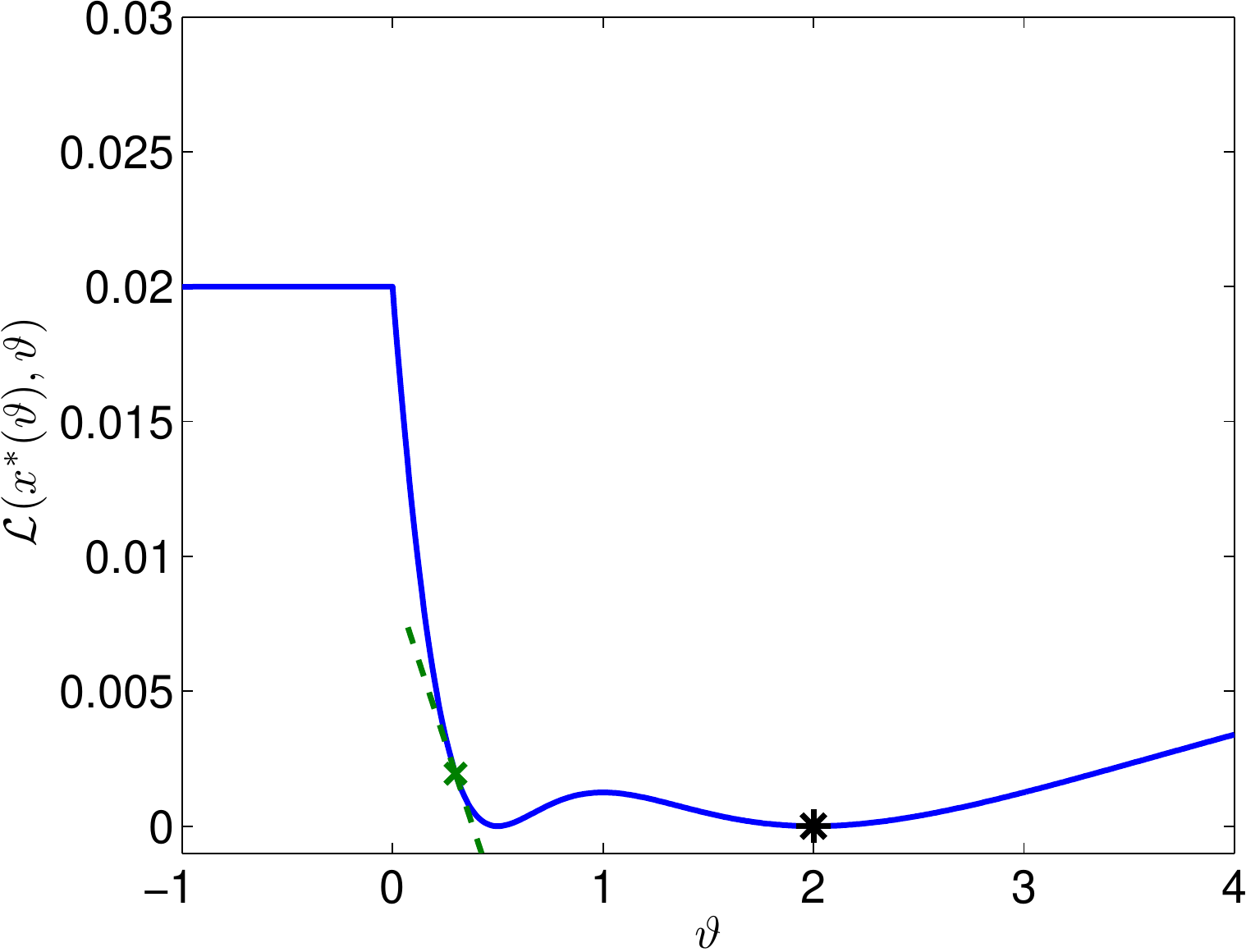}
\includegraphics[width=0.48\linewidth]{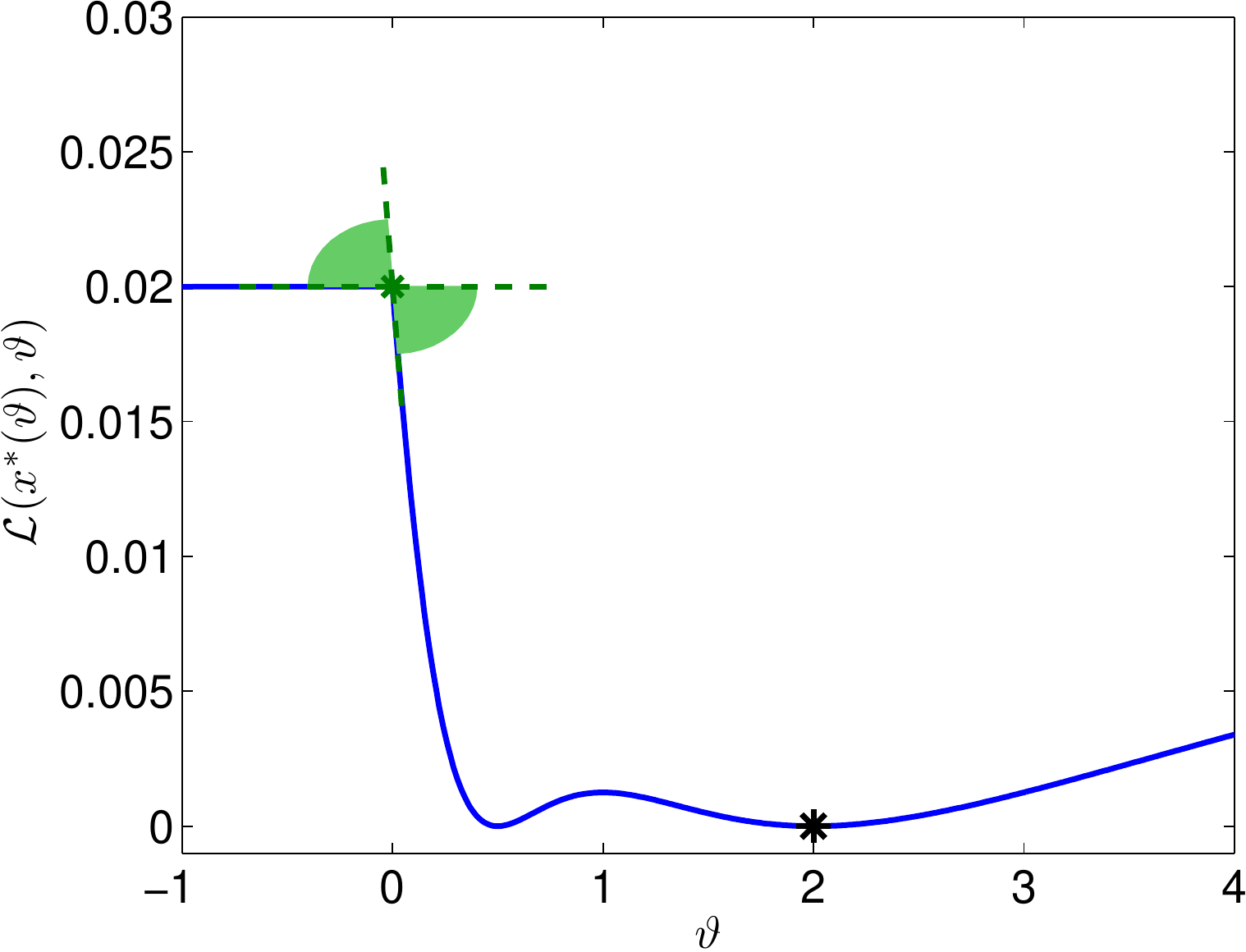}
\caption{\label{fig:example-nn-ls-tangents}Analytic tangents to the upper level objective function of \eqref{eq:bilevel-nonneg-LS} at $\theta=0.3$ and $\theta=0$. The function is non-smooth and, thus, at $\theta=0$ there exists many tangent lines.}
\end{figure}

In the experiments, we focus on the estimation of the gradient (in Figure~\ref{fig:example-nn-ls-tangents}). Therefore, the step size parameters of the individual algorithms are chosen such that a comparable convergence of the lower level energy is achieved, if possible. 

\begin{figure}[t]
\centering
\includegraphics[width=0.48\linewidth]{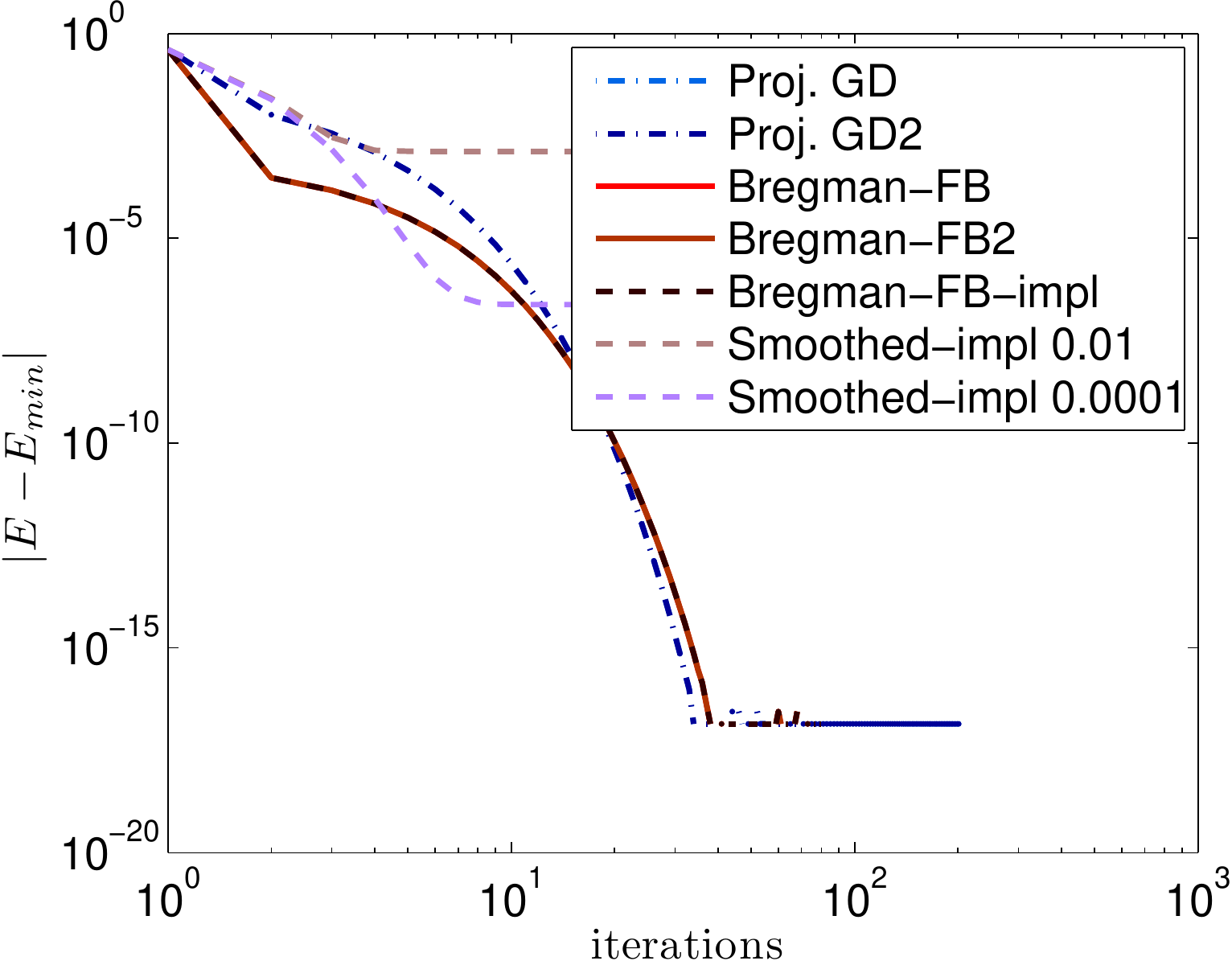}
\includegraphics[width=0.48\linewidth]{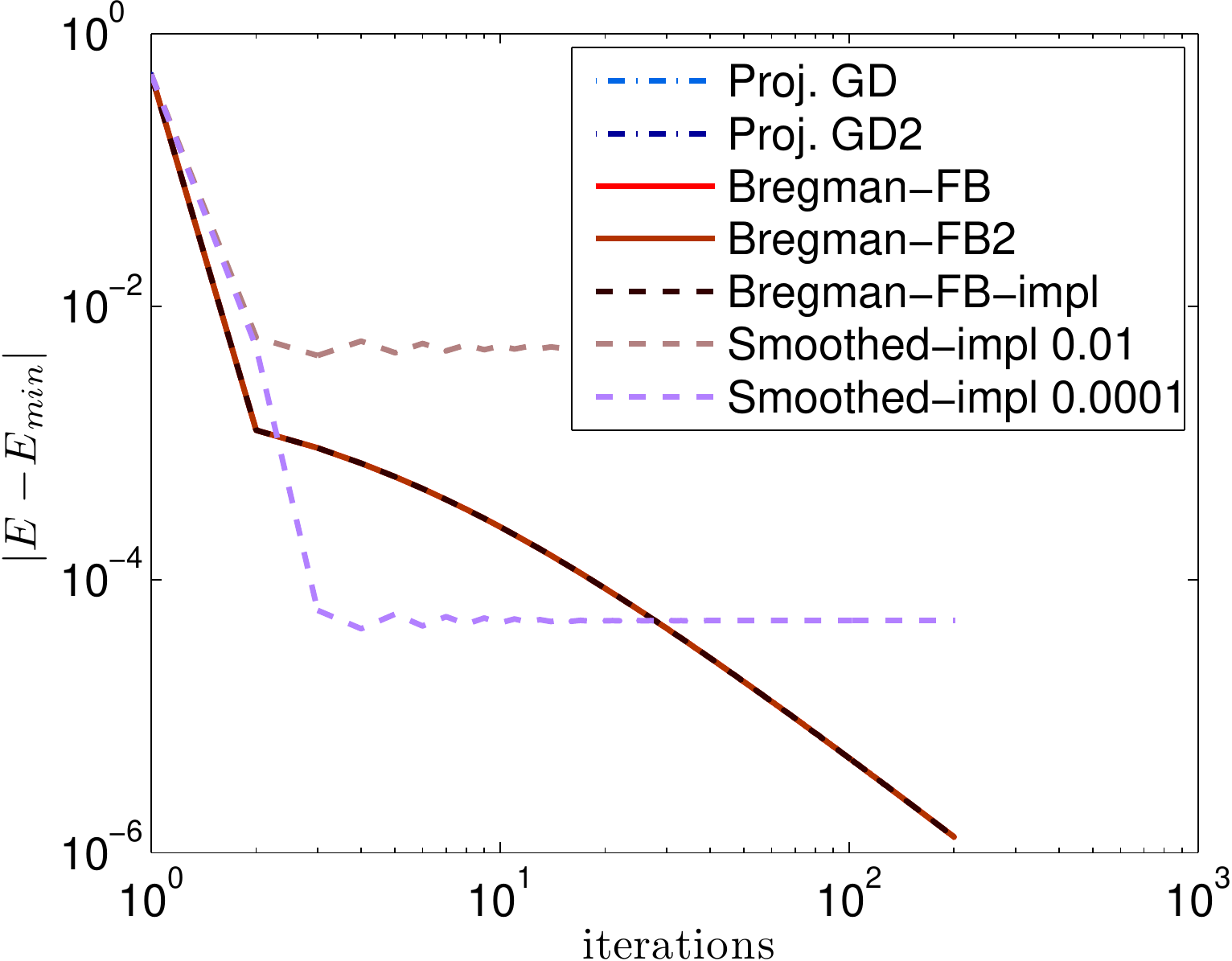}\\
\includegraphics[width=0.48\linewidth]{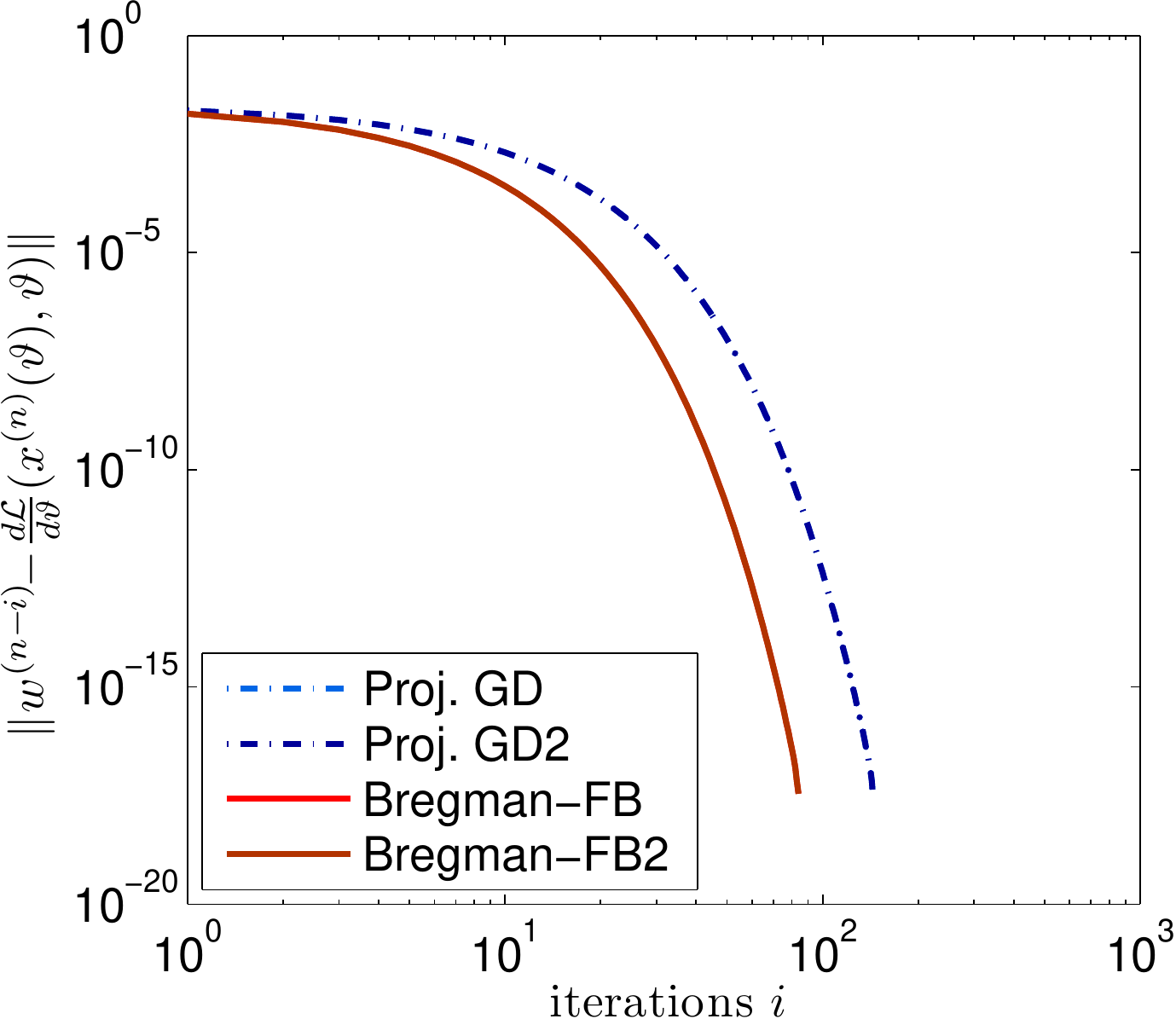}
\includegraphics[width=0.48\linewidth]{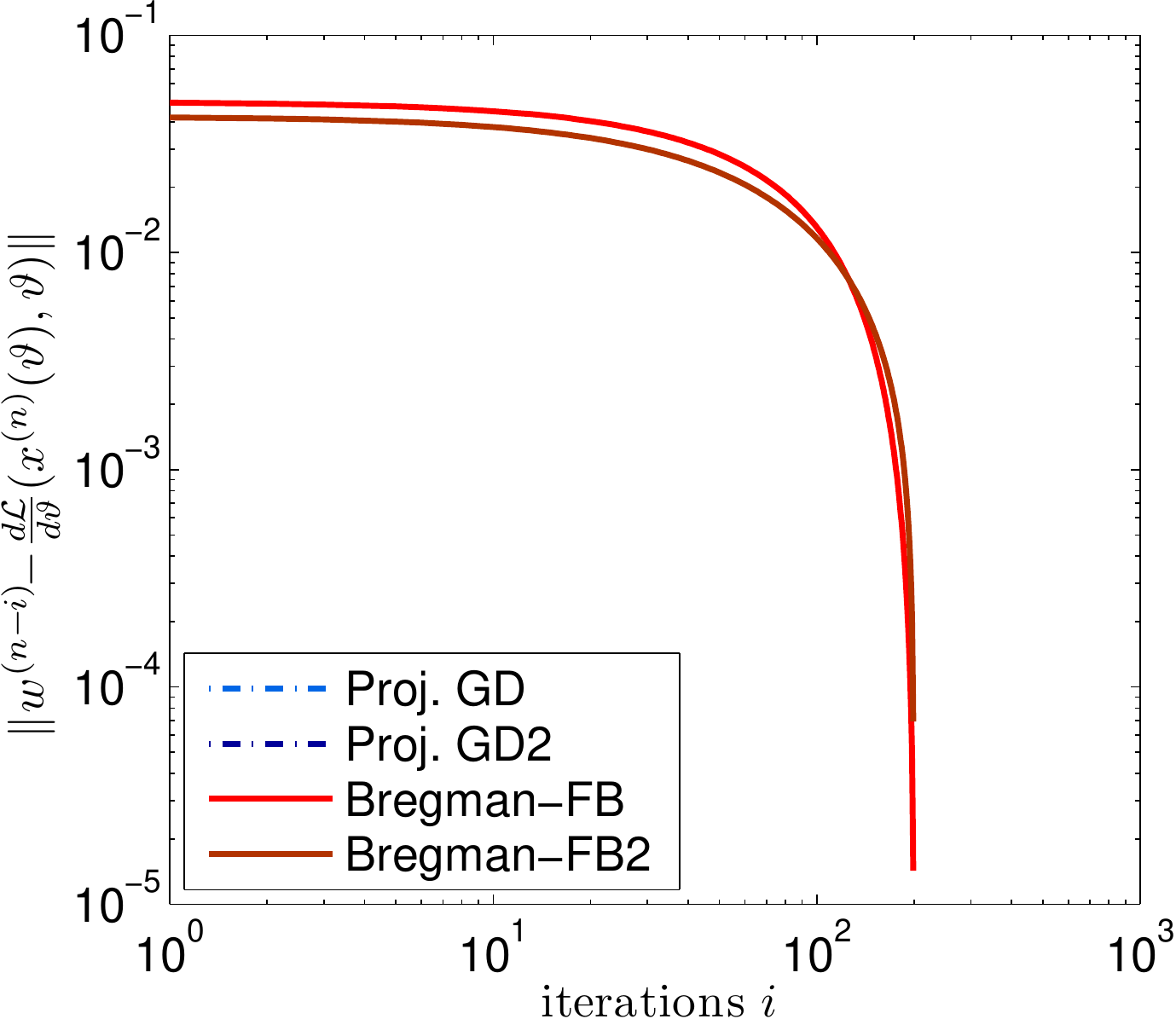}
\caption{\label{fig:example-nn-ls-grad-contrib}The upper row shows the energy decrease along the forward iterations. The lower row shows the convergence to the respective gradient value along the back-iterations. On the left hand side the plot is generated with $\theta=0.3$ and on the right hand side with $\theta=0$. The \enquote{-impl} methods do not appear in the bottom row as no back-iterations are involved. For $\theta=0$, due to the simple structure of the lower level problem, projected gradient descent converges exactly in one iteration, thus it is not shown. The gradient converges linearly to its final value, which means that often a few back-iterations are enough to achieve a gradient estimate of good quality. }
\end{figure}

For \ProjGD, \BregFB, and \BregFBb the chain rule must be applied recursively. We plot the change of the gradient accumulation along these back-iterations (of 200 forward-iterations) in bottom of Figure~\ref{fig:example-nn-ls-grad-contrib} and the energy evolution in the upper part of this figure. 
In this example, we can observe a linear decrease in the contribution to the respective final gradient value, which shows that back-iterations can be stopped after a few iterations without making large errors. 

Interestingly, the approximations \BregFBb and \ProjGDb work well, as they show the same gradient accumulation as \BregFB and \ProjGD, respectively. This situation changes when the number of forward-iterations is reduced. For about 15 forward-iterations, a difference of order $10^{-4}$ becomes visible (case $\theta=0.3$). 
\begin{figure*}[p]
\begin{center}
\includegraphics[width=0.32\linewidth]{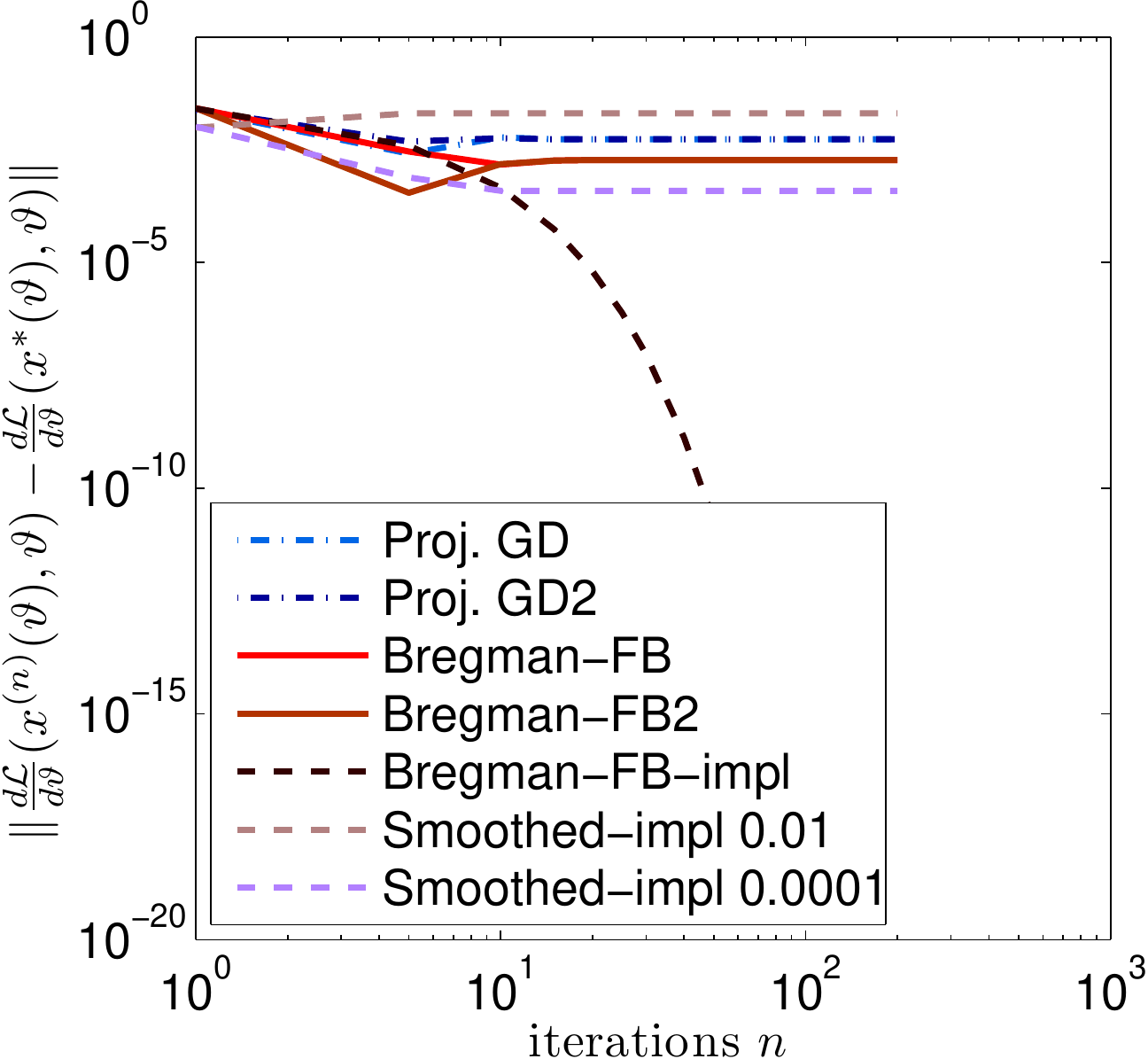}
\includegraphics[width=0.32\linewidth]{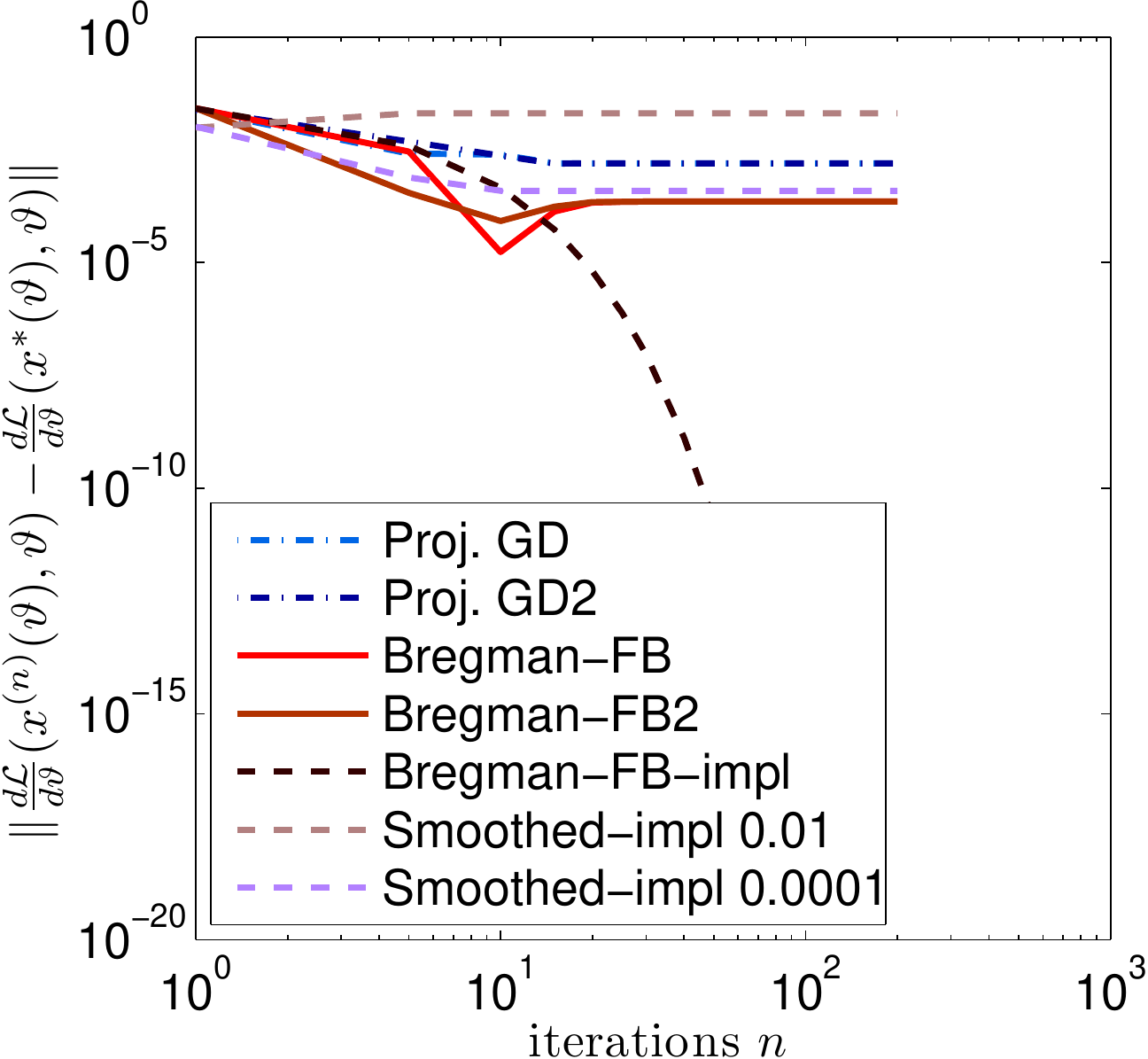}
\includegraphics[width=0.32\linewidth]{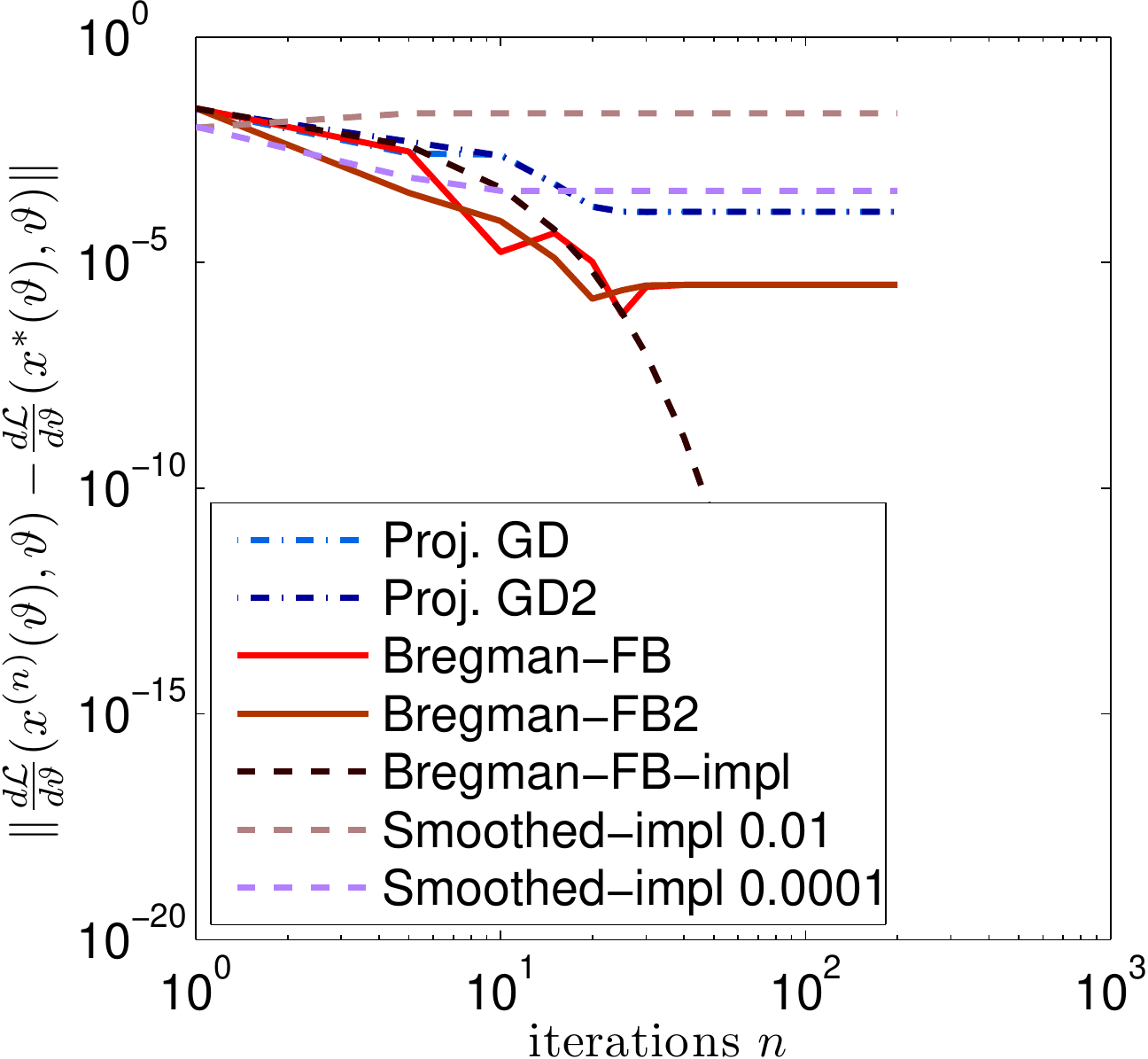}\\
\includegraphics[width=0.32\linewidth]{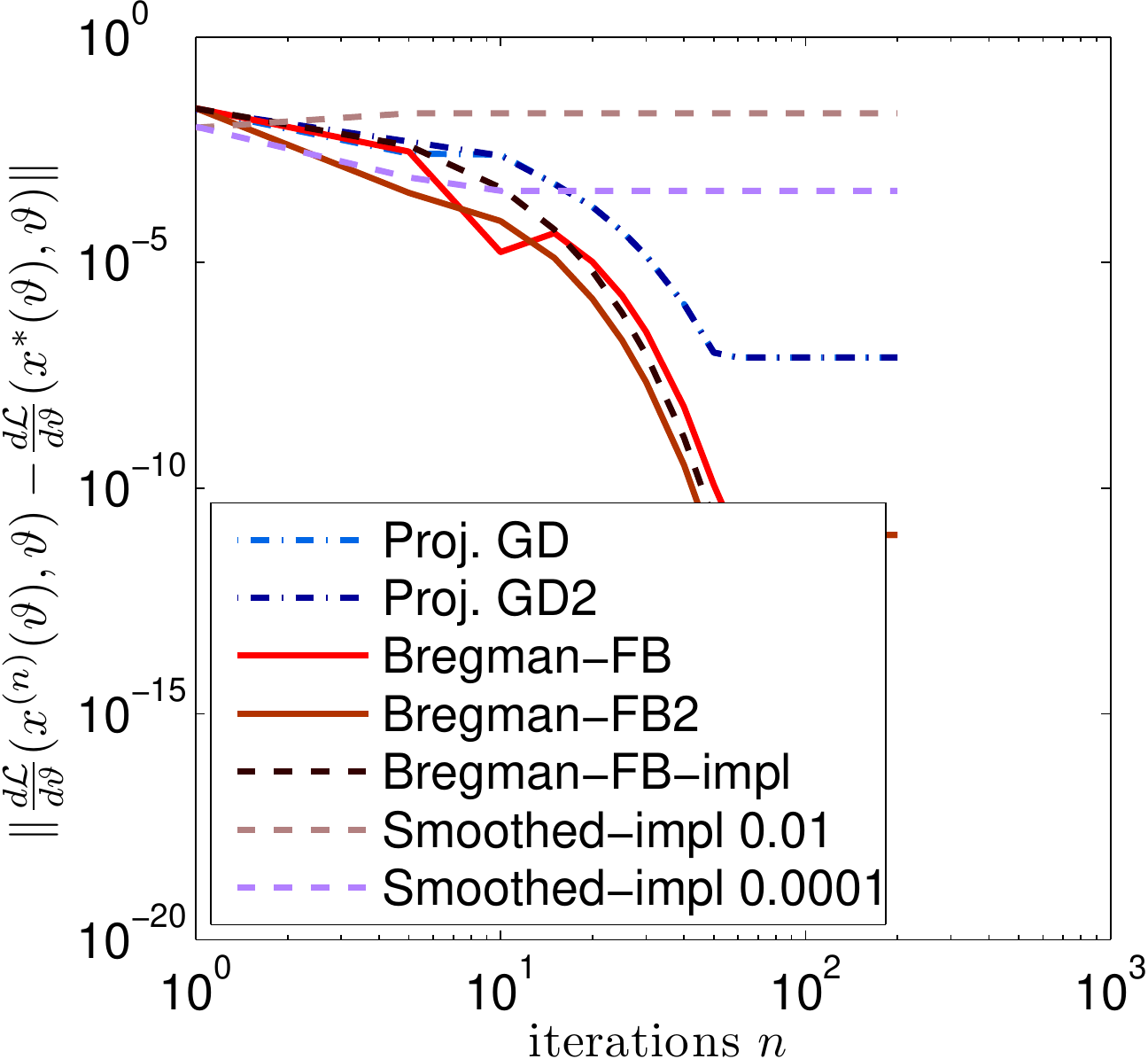}
\includegraphics[width=0.32\linewidth]{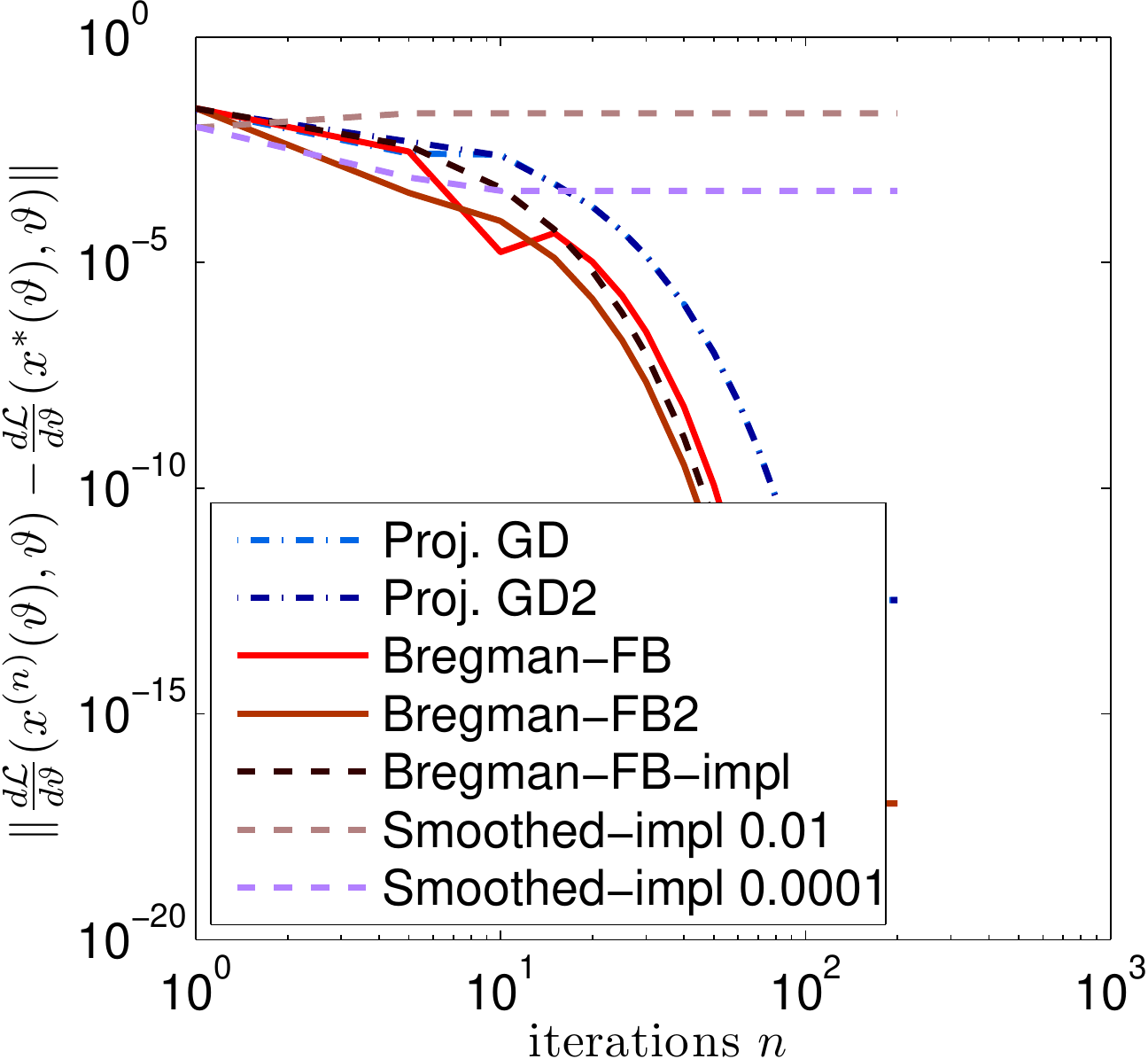}
\includegraphics[width=0.32\linewidth]{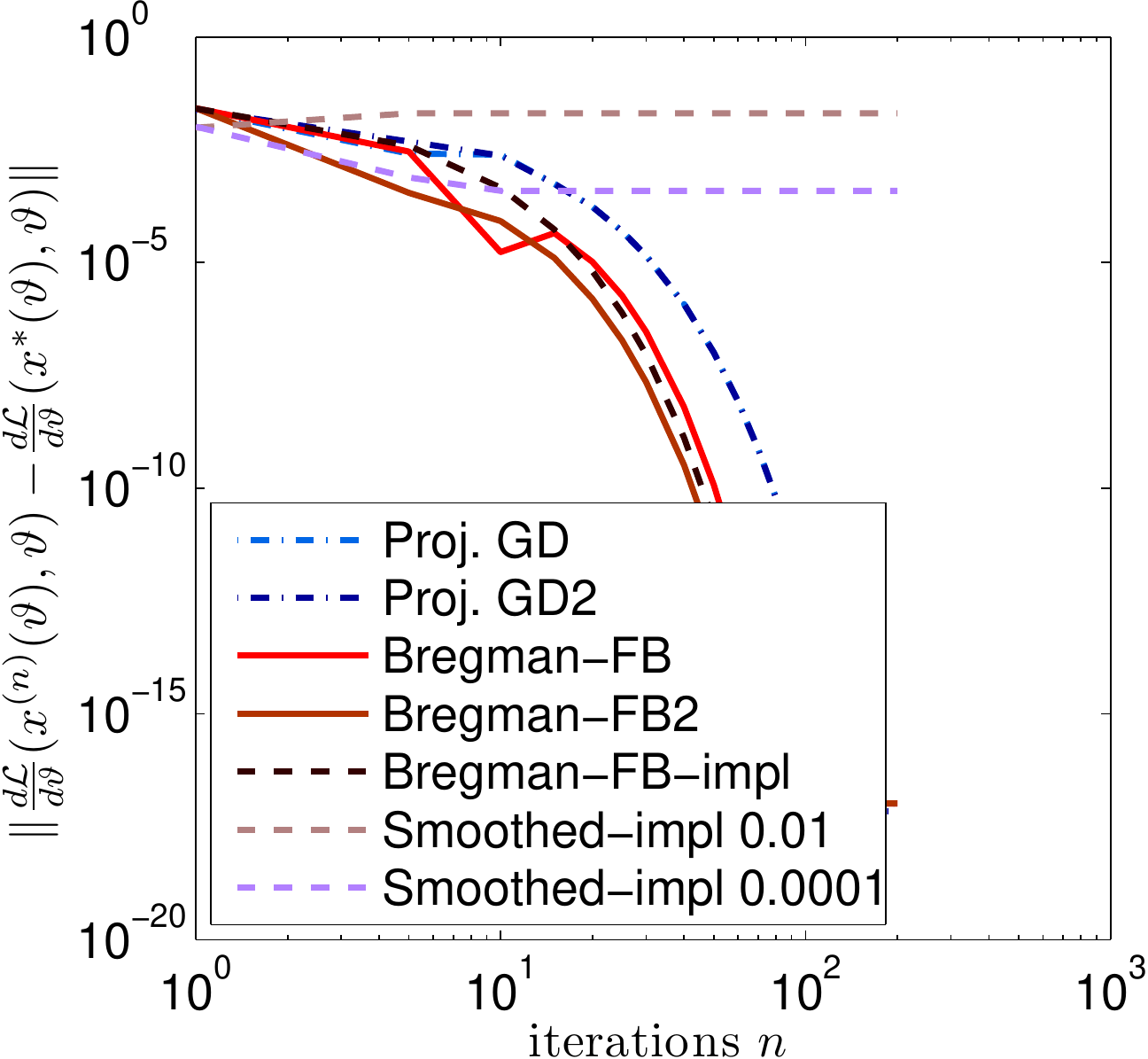}
\caption{\label{fig:example-nn-ls-grad-contrib-a}Convergence of the numerical gradients towards the analytic gradient for $\theta=0.3$. Row-wise, from left to right, the number of back-iterations is increased: 5, 10, 20, 50, 100, 200. More back-iterations lead to more accurate gradient estimates. The \enquote{-impl} methods always perform equally, as no back-iterations are required. \ImplDiff performs worst due to the rough approximation. Our methods \BregFB, \BregFBb, and \AlgImpl are the best; and converge slightly better than \ProjGD and \ProjGDb.}
\end{center}\medskip
\begin{center}
\includegraphics[width=0.32\linewidth]{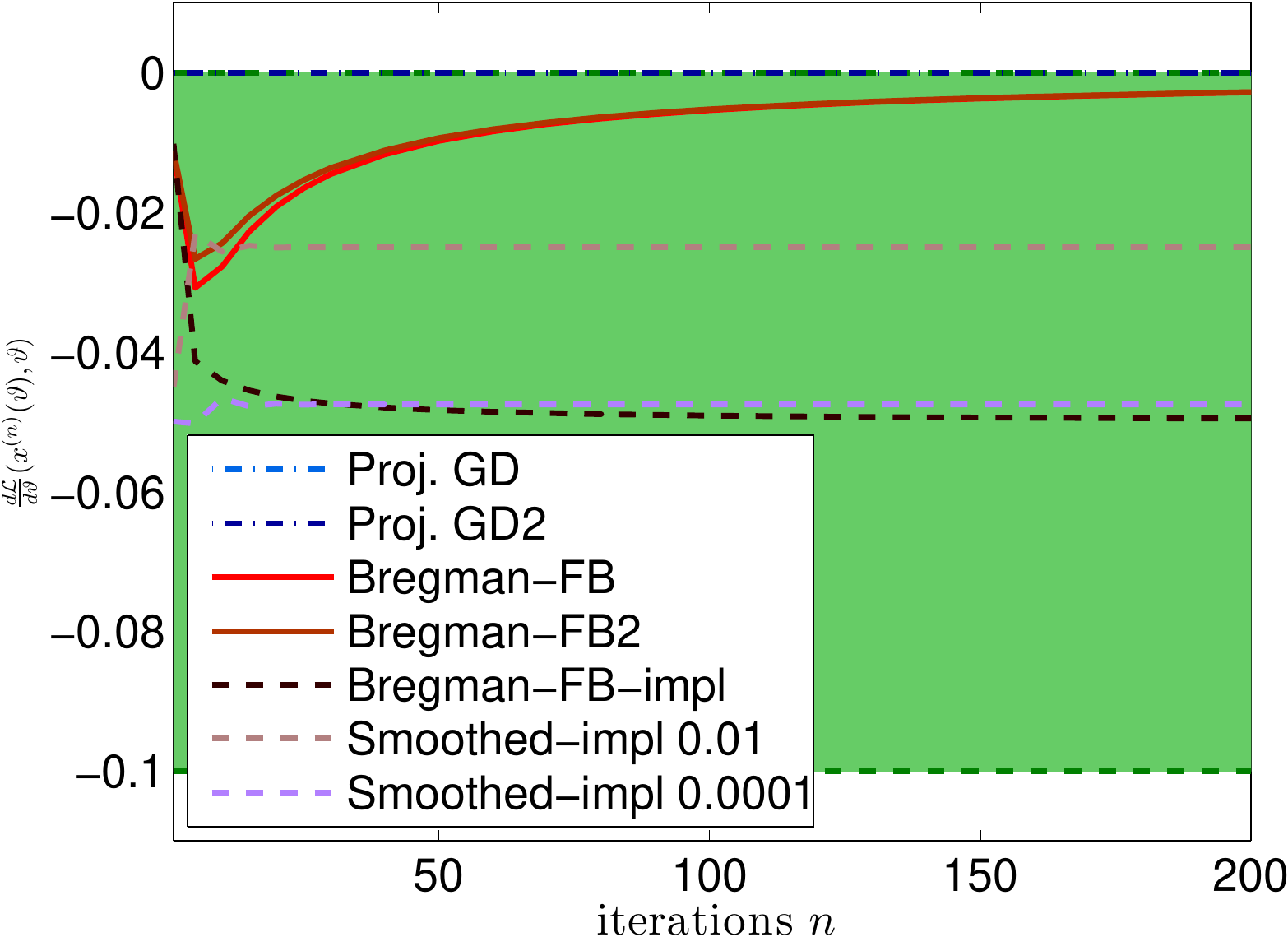}
\includegraphics[width=0.32\linewidth]{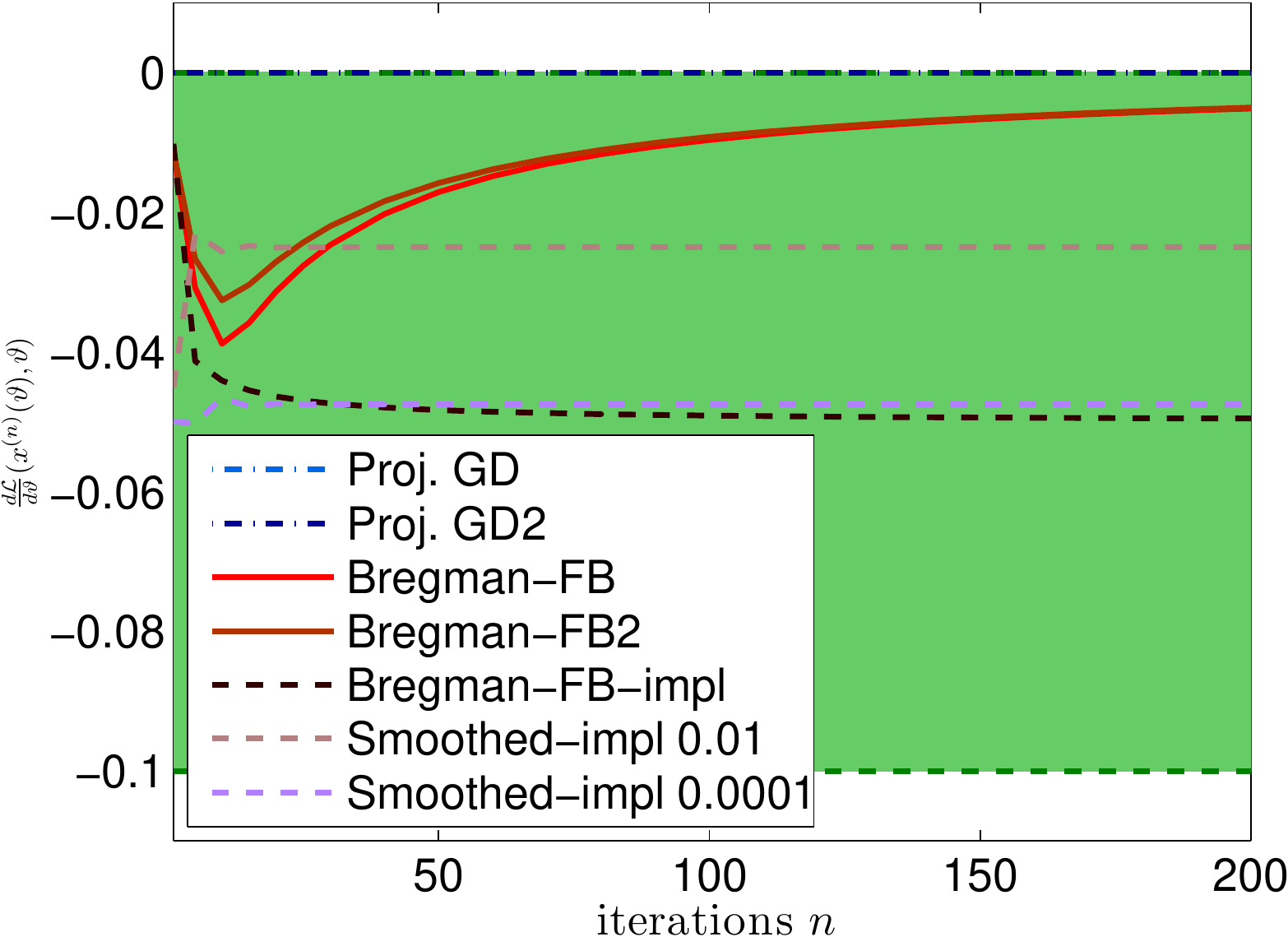}
\includegraphics[width=0.32\linewidth]{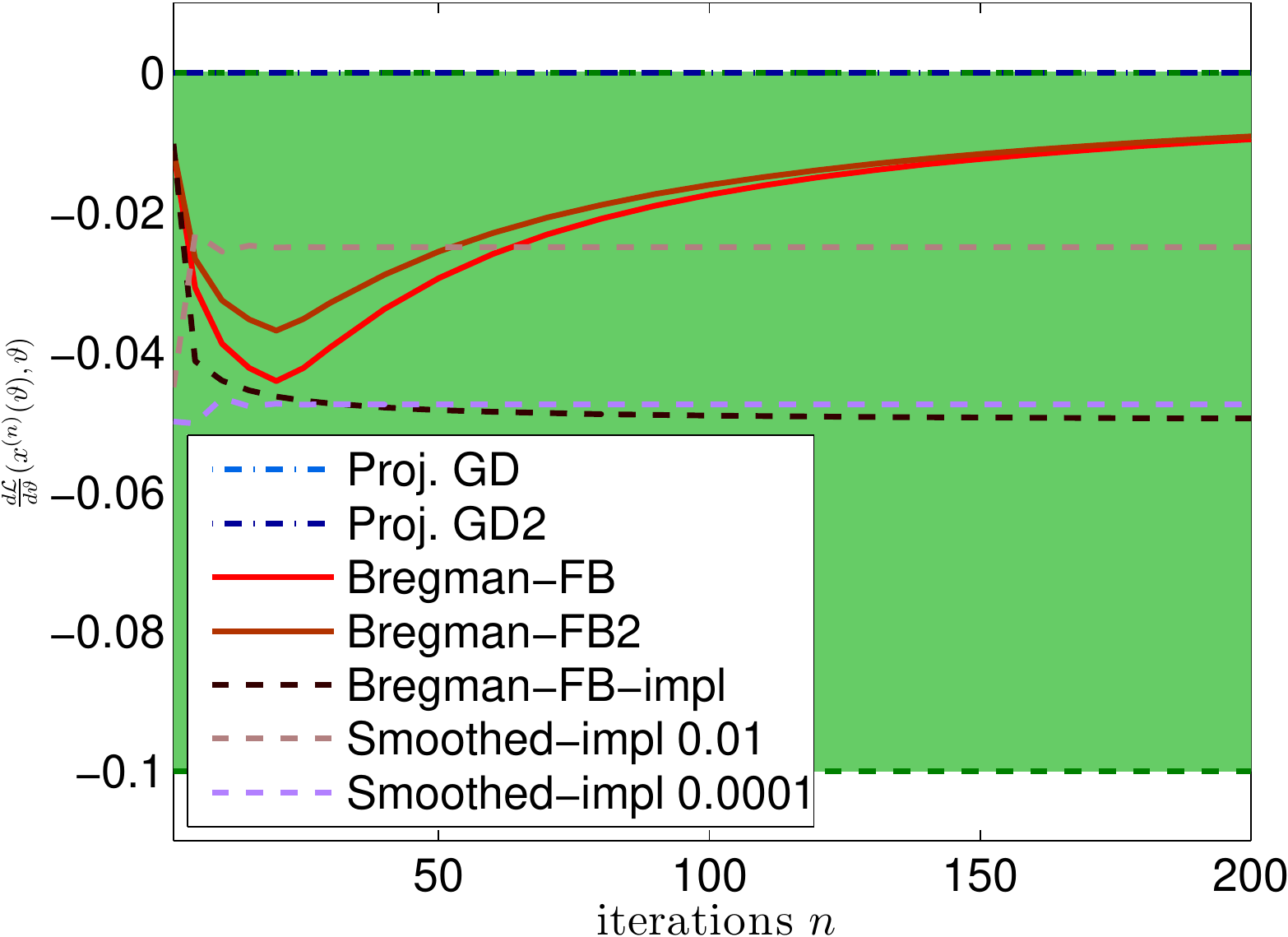}\\
\includegraphics[width=0.32\linewidth]{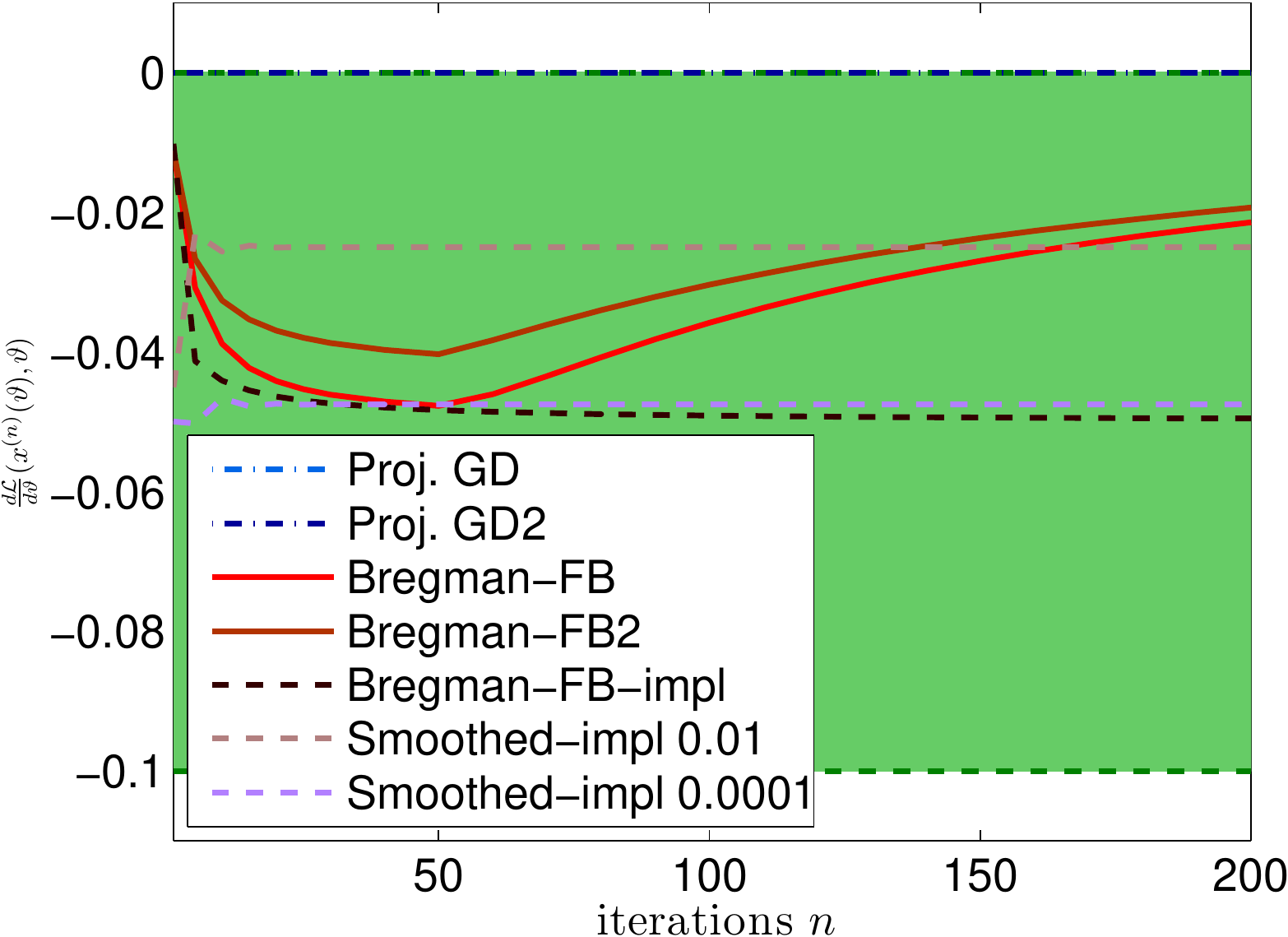}
\includegraphics[width=0.32\linewidth]{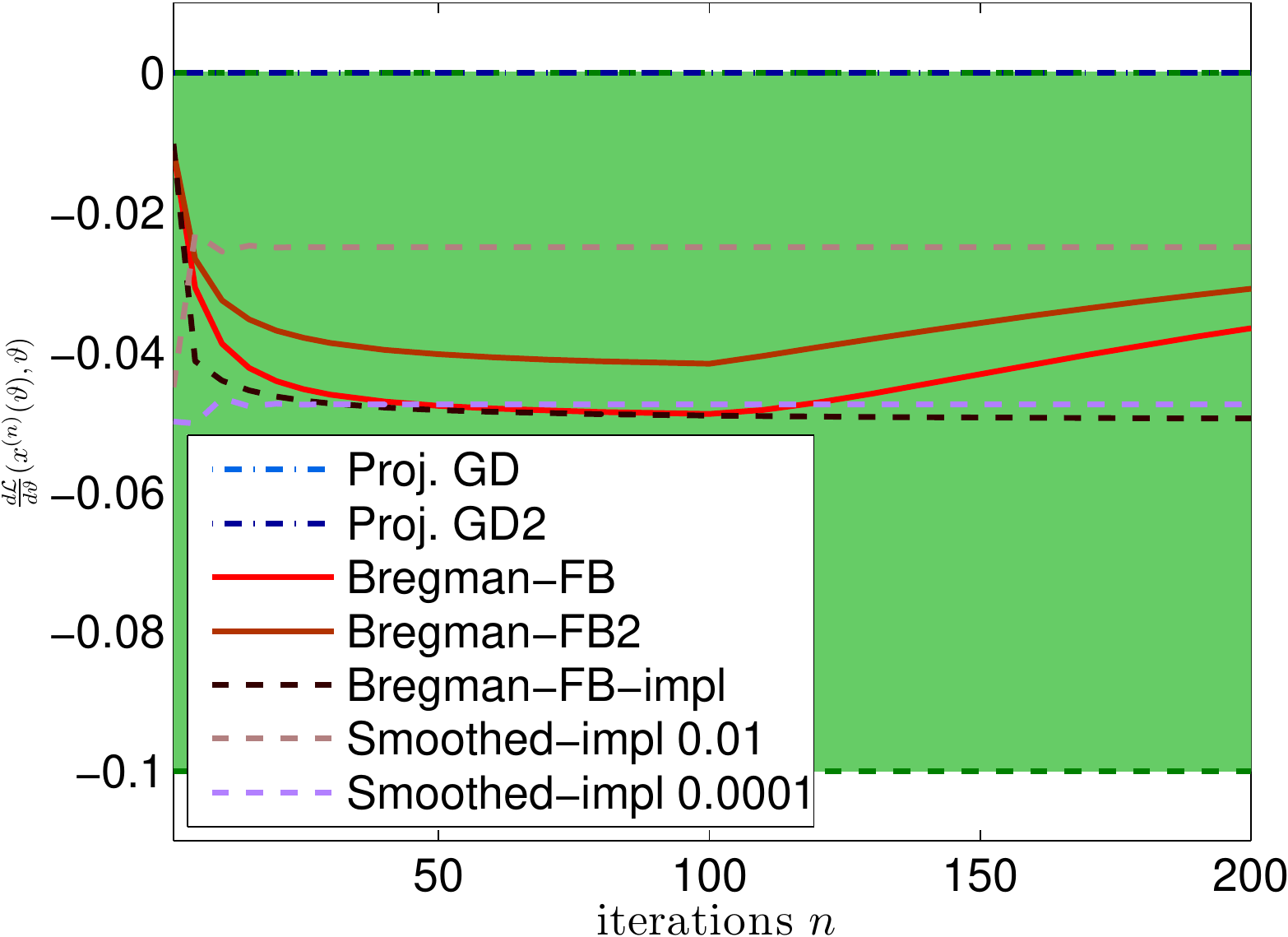}
\includegraphics[width=0.32\linewidth]{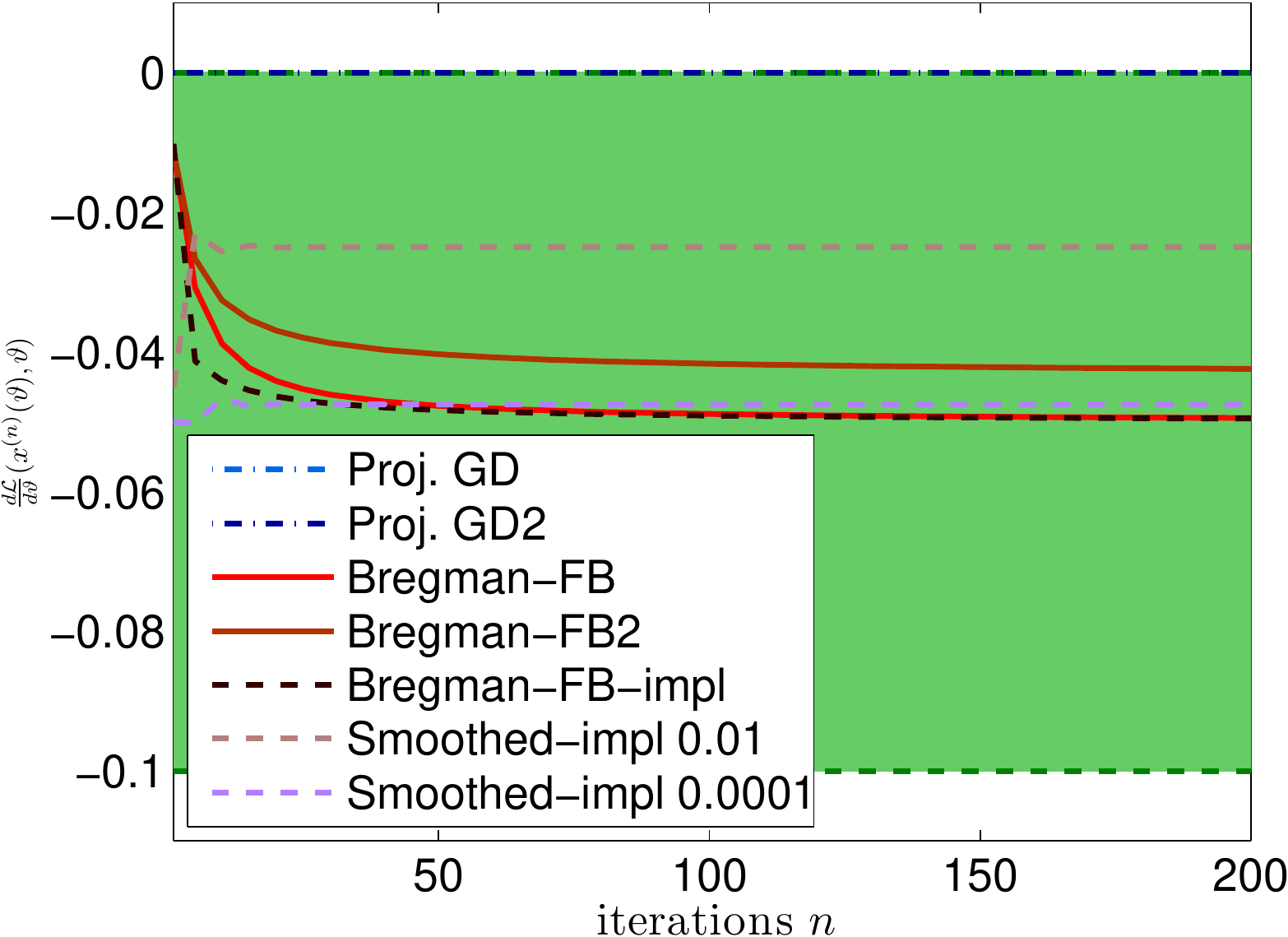}
\caption{\label{fig:example-nn-ls-grad-contrib-b}Convergence of the numerical gradients towards the analytic gradient for $\theta=0$. Row-wise, from left to right, the number of back-iterations is increased: 5, 10, 20, 50, 100, 200. All methods perform equally well, as they lie in the bright green area that indicates the range of the subdifferential.}
\end{center}
\end{figure*}

Figures ~\ref{fig:example-nn-ls-grad-contrib-a} and~\ref{fig:example-nn-ls-grad-contrib-b} address the convergence of the gradient towards the analytic gradient, with respect to different approximation accuracies (varied by the number of back-iterations). Figure~\ref{fig:example-nn-ls-grad-contrib-a} shows the convergence for $\theta=0.3$ and Figure~\ref{fig:example-nn-ls-grad-contrib-b} for $\theta=0$. Numerically, we observe convergence to the analytic gradients. 

Surprisingly, all methods perform equally well in the case $\theta=0$. The estimated gradient lies always in the subdifferential at this point. The range of the subdifferential is indicated with bright green color in Figure~\ref{fig:example-nn-ls-grad-contrib-b}. While \ProjGD and \ProjGDb estimate a gradient from the boundary of the subdifferential, the other methods estimate a subgradient from the interior. However, all of these values are feasible and belong to the analytic subdifferential. 

\section{Application to Multi-Label Segmentation} \label{sec:potts-seg-CNN}
In this section, we show how the idea can be applied in practice. To this end, we introduce a multi-label segmentation model. 
We use a convolutional neural network (CNN) to parametrize the segmentation model. Alternatively, this construction can be thought of as having a segmentation model as the final stage of a deep neural network.
In this setting, the bilevel problem amounts to finding the parameters of the CNN such that the loss on training data is minimized. The presented approach provides a generic way to train such systems in an end-to-end fashion.

\subsection{Model} \label{subsec:potts-model}
  Given a cost tensor $\cost\in\X^\nz$, where $\X=\R^{\nx\ny}$, that assigns to each pixel $(\px,\py)$ and each label $\pz$, $\px=1,\ldots,\nx$, $\py=1,\ldots,\ny$, $\pz=1,\ldots,\nz$, a cost $\cost^\pz_{\px,\py}$ for the pixel taking label $\pz$. We often identify $\R^{\nx\times\ny}$ with $\R^{\nx\ny}$ by $(\px,\py)\mapsto \px + (\py-1)\nx$ to simplify the notation. The sought segmentation $\optu\in \X_{[0,1]}^\nz$, where $\X_{[0,1]}=[0,1]^{\nx\ny}\subset\X$, is represented by a binary vector for each label. As a regularizer for a segment's plausibility we measure the boundary length using the total variation (TV). The discrete derivative operator $\map{\opK}{\X}{\Y}$, where we use the shorthand $\Y:=\X\times\X$ (elements from $\Y$ are considered as column vectors), is defined as:
  \[
  \begin{split}
    (\opK \optu^\pz)_{\px,\py} :=& \begin{pmatrix} (\opK \optu^\pz)_{\px,\py}^x \\ (\opK \optu^\pz)_{\px,\py}^y \end{pmatrix} \in \Y (= \R^{2\nx\ny}),\\
	\opKnz \optu :=& (\opK \optu^1,\ldots, \opK \optu^\nz), \\
    (\opK \optu^\pz)_{\px,\py}^x :=& \begin{cases}
        \optu^\pz_{\px+1,\py} - \optu^\pz_{\px,\py}\,, &\!\! \text{if } 1\leq \px < \nx, 1\leq \py\leq \ny \\
        0 \,, &\!\! \text{if } \px=\nx, 1\leq \py\leq \ny
    \end{cases} 
  \end{split}
  \]
  $(\opK \optu^\pz)_{\px,\py}^y$ is defined analogously. From now on, we work with the image as a vector indexed by $\ppi=1,\ldots,\nx\ny$. Let elements in $\Y$ be indexed with $\ppd=1,\ldots,2\nx\ny$. Let the inner product in $\X$ and $\Y$ be given, for $u^\pz,v^\pz\in\X$ and $p^\pz,q^\pz\in\Y$, as:
  \[
  \begin{split}
    \scal{u^\pz}{v^\pz}_\X &:= \sum_{\ppi=1}^{\nx\ny} u^\pz_{\ppi} v^\pz_{\ppi}, \;
    \scal{p^\pz}{q^\pz}_\Y := \sum_{\ppd=1}^{2\nx\ny} p^\pz_{\ppd} q^\pz_{\ppd},\\
  	\scal{u}{v}_{\X^\nz} &:= \sum_{\pz=1}^\nz \scal{u^\pz}{v^\pz}_\X, \; 	\scal{p}{q}_{\Y^\nz} := \sum_{\pz=1}^\nz \scal{p^\pz}{q^\pz}_\Y.
  \end{split}
  \] 
  The (discrete, anisotropic) TV norm is given by 
  \[
    \norm[1]{\opKnz \optu} := \sum_{\pz=1}^\nz \sum_{\ppd=1}^{2\nx\ny} \abs{ (\opK \optu^\pz)_{\ppd} }\,,
  \]
  where $\abs\cdot$ is the absolute value. In the following, the variables $\ppi=1,\ldots,\nx\ny$ and $\ppd=1,\ldots,2\nx\ny$ always run over these index sets, thus we drop the specification; we adopt the the same convention for $\pz=1,\ldots,\nz$. We define the pixel-wise nonnegative unit simplex
  \begin{align} \label{eq:unit-simplex}
    \Simplex^\nz := \{\forall & (\ppi,\pz)\colon 0\leq \optu^\pz_{\ppi}\leq 1 \nonumber\\ &\text{ and } \forall \ppi\colon \textstyle\sum_\pz \optu^\pz_\ppi = 1 \;{\optu\in\X^\nz}\}\,,
  \end{align}
  and the pixel-wise (closed) $\ell_\infty$-unit ball around the origin
  \[
    B^{\ell_\infty}_{1}(0) := \set[{\forall (\ppd,\pz)\colon \abs{\optp^\pz_{\ppd}} \leq 1}]{\optp\in \Y^\nz} \,.
  \]
  Finally, the segmentation model reads
  \begin{equation} \label{eq:seg-potts}
    \min_{\optu\in\X^\nz}\ \scal{\cost}{\optu}_{\X^\nz} + \norm[1]{\WW \opKnz \optu}\,,\quad \st\ \optu \in\Simplex^\nz \,,
  \end{equation}
  where we use a diagonal matrix $\WW$ to support contrast-sensitive penalization of the boundary length.
  
  This model and the following reformulation as a saddle-point problem are well known (see e.g. \cite{CP11})
  \begin{align} \label{eq:seg-potts-saddle-point}
    \min_{\optu\in\X^\nz} & \max_{\optp\in\Y^\nz}\ \scal{\WW \opKnz \optu}{\optp}_{\Y^\nz} + \scal{\optu}{\cost}_{\X^\nz} \,, \\
   & \st\ \optu\in\Simplex^\nz,\ \optp \in B^{\ell_\infty}_1(0)\,. \nonumber
  \end{align}
  The saddle-point problem \eqref{eq:seg-potts-saddle-point} can be solved using the ergodic primal-dual algorithm \cite{CP15}, which leads to an iterative algorithm with totally differentiable iterations. The primal update in \eqref{eq:PD-abstract} is discussed in Example~\ref{ex:Bregman-prox-lin-simplex} and the dual update of \eqref{eq:PD-abstract} is essentially Example~\ref{ex:Bregman-prox-lin-interval}. As a consequence Algorithm~\ref{alg:bilevel-pd-alg-reverse} can be applied to estimate the derivatives. A detailed derivation of the individual steps of the algorithm can be found in \cite{ORBP15}.
 
 \begin{figure*}[ht]
 \centering
 \includegraphics[width=0.319\linewidth]{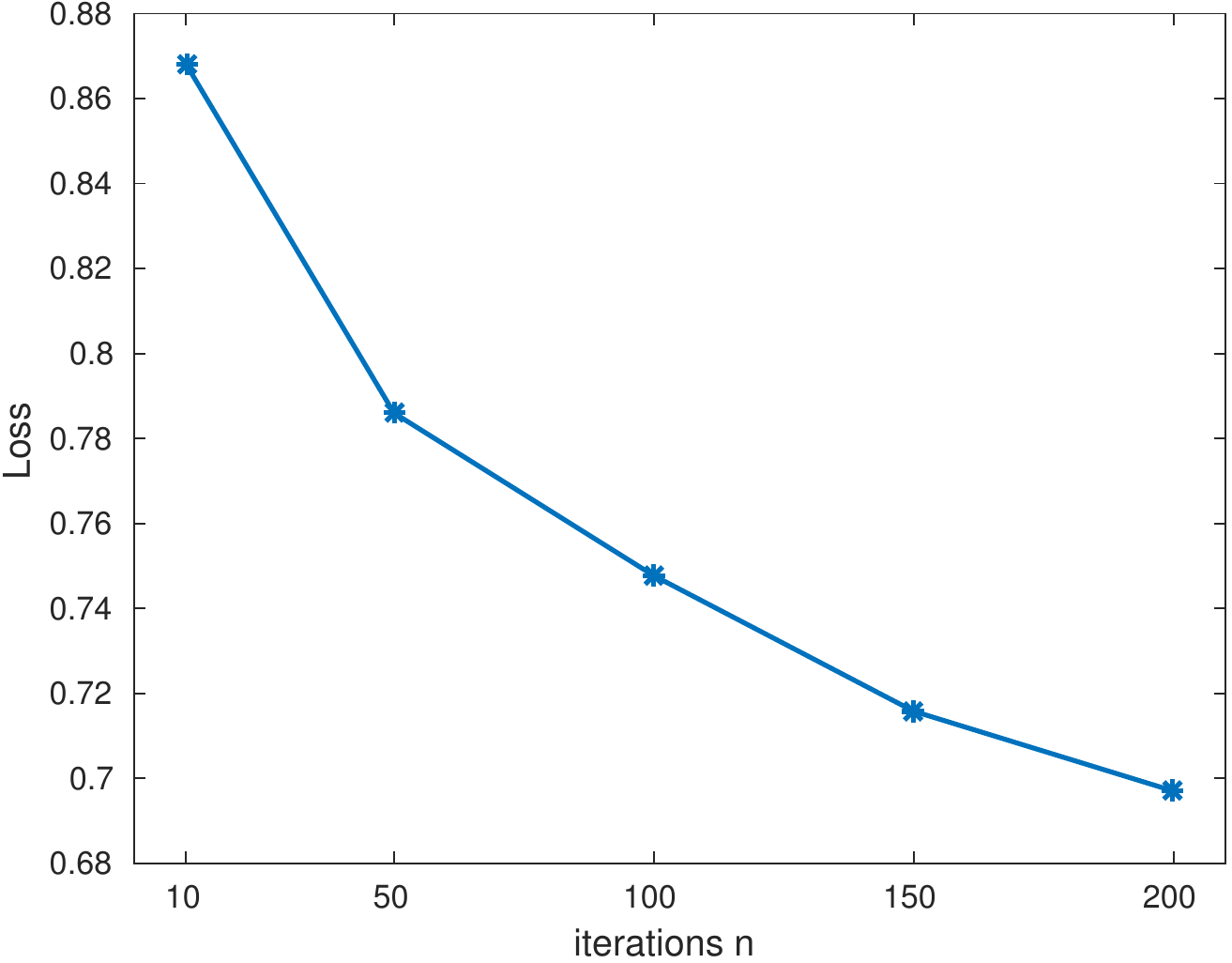}~
 \includegraphics[width=0.31\linewidth]{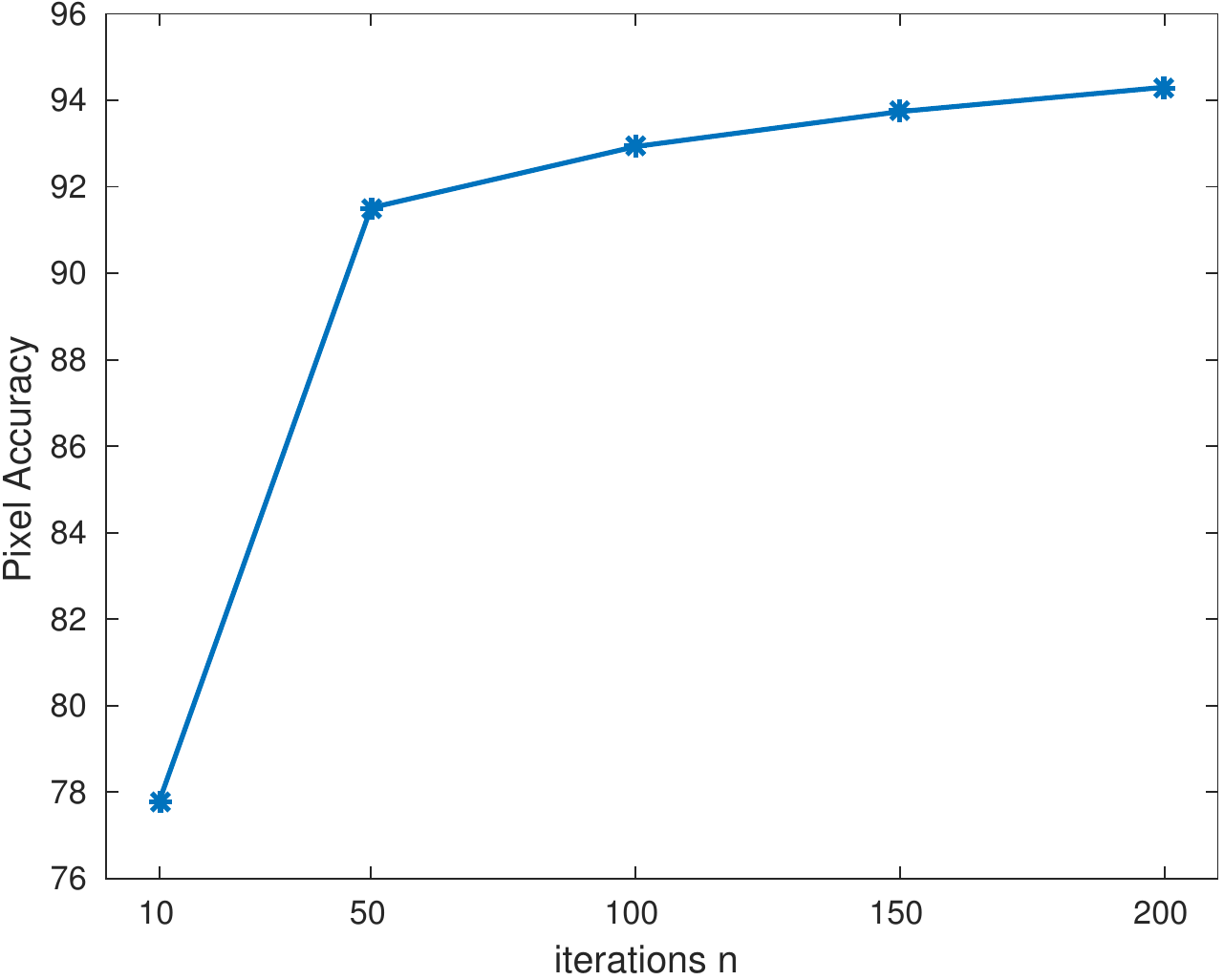}~
 \includegraphics[width=0.31\linewidth]{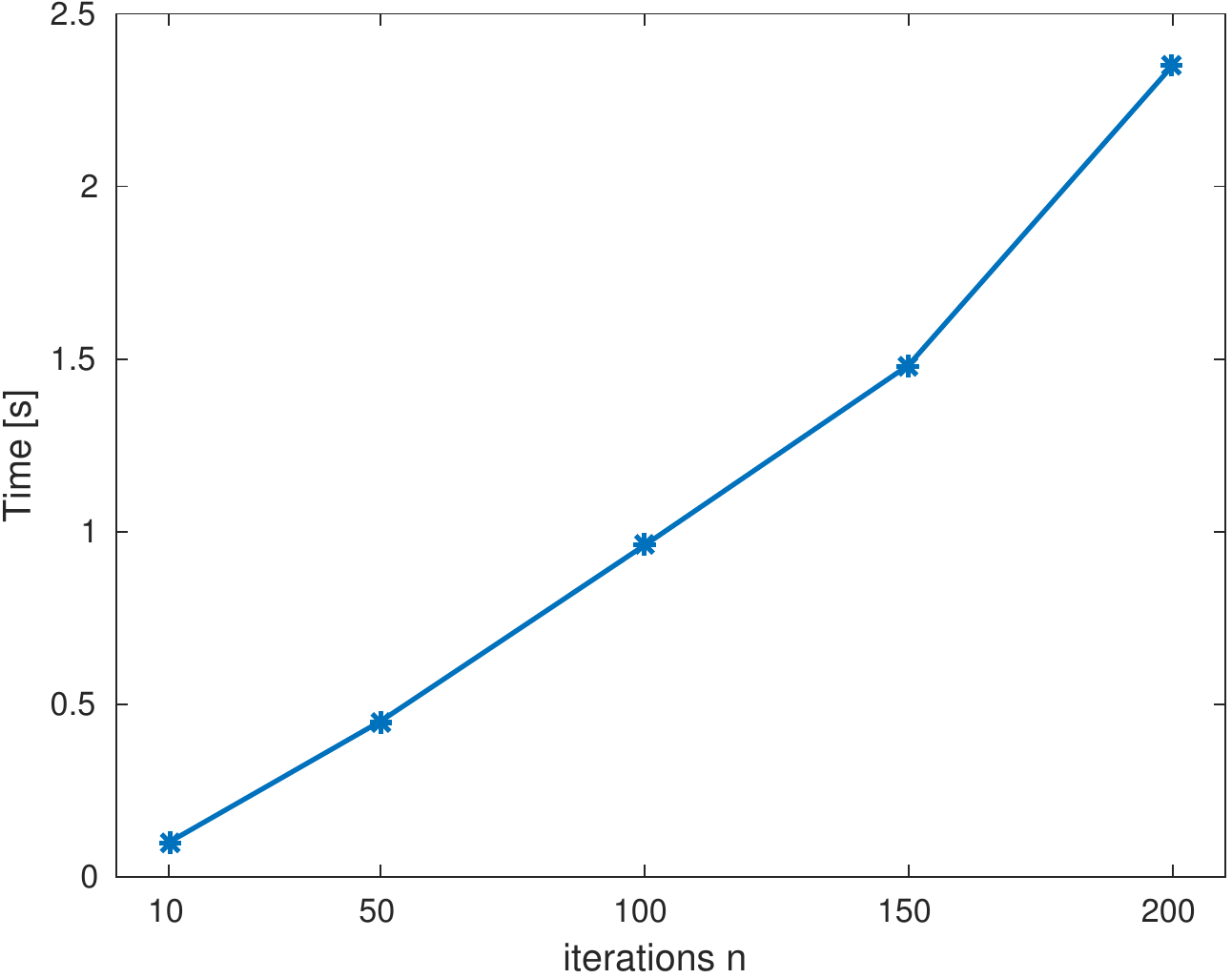}
 \caption{\label{fig:multilabel_iters} Training error vs. number of iterations of the algorithm solving the lower level problem. From left to right, average per-pixel loss, per-pixel accuracy and time per outer iteration. The timing includes the forward pass as well as the gradient computations. Timings were taken on a NVIDIA Geforce Titan X GPU. A higher number of iterations clearly leads to lower error, but comes at the cost of a higher computational complexity.}
 \end{figure*}

\subsection{Parameter Learning} \label{subsec:param-learn-binary-seg}
We consider \eqref{eq:seg-potts} where the cost $\cost$ is given by the output of a CNN which takes as input an image $\img\in X^\nc$ to be segmented and is defined via a set of weights $\theta$. Formally, we have $\cost^\pz_{\ppi} = \cnn^\pz_{\ppi}(\theta,\img)$ with $\cnn :  \R^{N_\theta} \times X^\nc \rightarrow X^\nz$, where $\nc$ denotes the number of channels of the input image and $N_\theta$ is the number of weights parametrizing the CNN.

The training set consists of $\nt$ images
$\img^{1},\ldots,\img^{\nt}\in \X^\nc$ and their corresponding ground truth
segmentations $\gt^1,\ldots,\gt^{\nt} \in \{1,\ldots,\nz\}^{\nx\ny}$.

The parameters $\theta$ of the CNN are cast as an instance of the general bilevel optimization problem \eqref{eq:bilevel-general}:
  \begin{align} \label{eq:model-bilevel-seg}
     \min_{\theta\in\R^{N_\theta}}\ &
	\sum_{\pt=1}^\nt  \sum_{\ppi=1}^{\nx\ny}\log \Big(\sum_{\pz=1}^\nz \exp(\optu^\pz_\ppi(\theta,\img^\pt))\Big) %
  - \gt^\pt_\ppi(\theta,\img^\pt)\nonumber\\
     & \st\ \optu(\theta,\img^\pt) = \arg\min_{\optu\in\X^\nz} E(\optu, \cnn(\theta,\img^\pt)),
  \end{align}
where energy $E$ in the lower level problem is \eqref{eq:seg-potts} and the higher-level problem is defined as the softmax loss.
\begin{remark}
We could equivalently use a multinomial logistic loss, since $\optu_\ppi(\theta, \img^\pt)$ lies in the unit simplex by construction.  We use this definition to allow for a simplified treatment of the case of training a CNN without the global segmentation model.
\end{remark}

\subsection{Experiments}

We implemented our approach as a custom layer in the MatConvNet framework \cite{matconvnet}. We used the Stanford Background dataset \cite{stanford_background}, which consists of 715 input images and pixel-accurate ground truth consisting of the geometric classes \textit{sky}, \textit{vertical} and \textit{horizontal}. We used ADAM \cite{adam} for the minimization of the higher-level problem. We found that general plain stochastic gradient descent performs poorly in our setting, since the global segmentation model can lead to vanishing gradients.

In a first experiment we used a small subset of 9 images from the dataset to show the influence of the number of iterations used to solve the lower-level problem \eqref{eq:seg-potts} on the training objective. 
We learned a small network consisting of four layers of alternating convolutions with a kernel width of 3 pixels and ReLU units followed by a fully connected layer. We added $3\times 3$ max-pooling layers with a stride of two after the first and the second convolutional layers, which effectively downsamples the responses by a factor of 4. We added an up-convolutional layer to upsample the responses to the original image size. The penultimate layer of the CNN consist of a multiplicative scaling (analogous to a scalar smoothness parameter) of the CNN output followed by the global segmentation model \eqref{eq:seg-potts}.
We ran ADAM with a learning rate of $10^{-3}$ for a total of 1000 iterations with a mini-batch size of one image to learn the parameters of this network.

\begin{figure*}[ht!]
\centering
\includegraphics[width=0.245\linewidth]{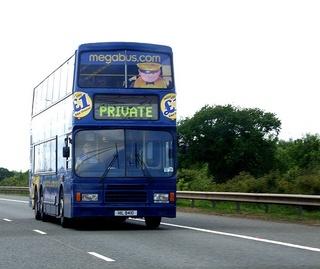}
\includegraphics[width=0.245\linewidth]{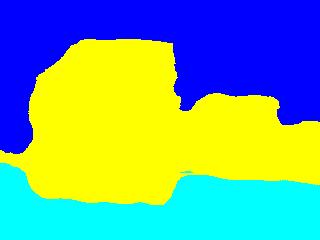}
\includegraphics[width=0.245\linewidth]{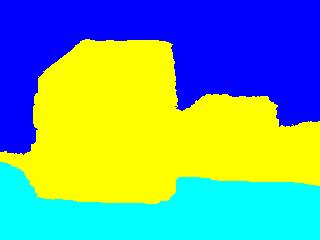}
\includegraphics[width=0.245\linewidth]{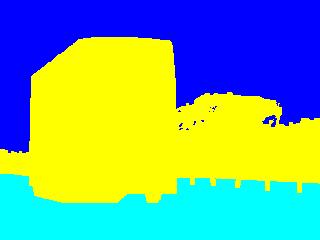}\\
\vspace{0.1cm}
\includegraphics[width=0.245\linewidth]{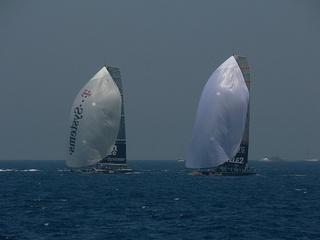}
\includegraphics[width=0.245\linewidth]{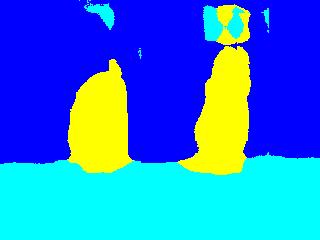}
\includegraphics[width=0.245\linewidth]{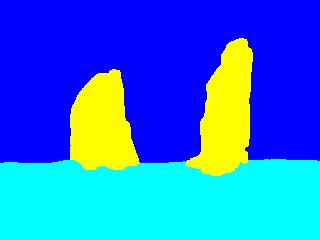}
\includegraphics[width=0.245\linewidth]{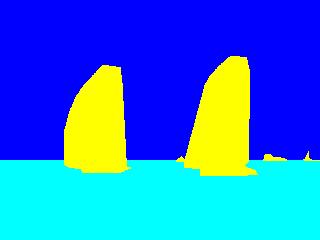}
\caption{\label{fig:multilabel_img} Example results from the test set. Row-wise, from left to right: Input image, CNN, CNN+Global, ground truth. The global model is able to align results to edges and is able to correct spurious errors.}
\label{fig:stanfordresimg}
\end{figure*}

Figure \ref{fig:multilabel_iters} shows the average per-pixel loss, the average pixel accuracy as well as the time per ADAM iteration  vs. number of iterations used to solve the lower-level problem (inner iterations).
This experiment shows that by solving the lower-level problem to higher accuracy the overall capacity and thus the accuracy of the system can be enhanced. This comes at the price of a higher computational complexity, which increases linearly with the number of iterations.

Finally, we performed a large scale experiment on this dataset. We partitioned the images into a training set of 572 images and used the remaining 143 images for testing. We used the pre-trained  Fully Convolutional Network \textit{FCN-32s}~\cite{long_shelhamer_fcn} as basis for this experiment. We adapted the geometry of the last two layers to this dataset and retrained the network. We then added a multiplicative scaling layer followed by the global segmentation model and refined the parameters. The number of inner iterations was set to 100, which provides a good trade-off between accuracy and computational complexity. We use a mini-batch size of 5 images and a learning rate of $10^{-3}$.

The average accuracy in terms of the average pixel accuracy (\textit{Acc}) in percent and Intersection over Union (\textit{IoU}) on both the test and the training set is shown in Table \ref{tab:stanfordres}. We compare the plain Fully Convolutional Network \textit{FCN} to the network with the additional global segmentation model \textit{FCN+Global}. We observed an increase of $1.4 \%$ in terms of IoU on the test set when using the global model.  This can be attributed to the fact that the CNN alone already provides good but coarse segmentations and the segmentation model uses only simple pairwise interactions. As such it is unable to correct gross errors of the CNN.

Since the presented approach is applicable to a broad range of energies,
training of more expressive energies which include more complex interactions (cf. \cite{crfasrnn}) is a promising direction of future research. Example segmentations from the test set are shown in Figure~\ref{fig:stanfordresimg}.

\begin{table}[b]
  \centering
  \begin{tabular}{lcccc}
  \toprule 
  & \multicolumn{2}{c}{Test} & \multicolumn{2}{c}{Train} \\
  & Acc & IoU & Acc & IoU \\ 
  \midrule FCN& 92.40 & 82.65 & 97.54  & 92.21 \\ 
  FCN+Global & \textbf{93.00} & \textbf{84.01} & \textbf{97.90}  & \textbf{93.53} \\ 
  \bottomrule 
  \end{tabular}
 \caption{\label{tab:stanfordres}Accuracy on the Stanford Background dataset \cite{stanford_background}. We compare the plain CNN to the CNN with an additional global segmentation model.}
\end{table}

\begin{remark}
  For a comparison to the smoothing approach from Section~\ref{sec:derivative-impl-fun} we refer to the conference version \cite{ORBP15}. 
\end{remark}

\section{Conclusion} \label{sec:conclustion}

We considered a class of bilevel optimization problems with non-smooth lower level problem. By an appropriate approximation we can formulate an algorithm with a smooth update mapping that solves a non-smooth optimization problem in the lower level. This allows us to apply gradient based methods for solving the bilevel optimization problem. 
A second approach directly considers the fixed-point equation of the algorithm as optimality condition for the lower level problem. Key for both ideas are Bregman proximity functions.

The idea of estimating gradients for an abstract algorithm was exemplified for a forward--backward splitting method and a primal--dual algorithm with Bregman proximity functions. Several potential application examples were shown. A toy example confirmed our results and provided some more intuition. The contribution of our idea to practical applications was demonstrated by a multi-label segmentation model that was coupled with a convolutional neural network.

There are several open questions, for example convergence of the sequence of gradients or a full classification of optimization problems that allow for algorithms with smooth update mapping.


\subsection*{Acknowledgment}
Peter Ochs and Thomas Brox acknowledge support from the German Research
Foundation (DFG grant BR 3815/8-1). Ren\'{e} Ranftl acknowledges support from Intel Labs. Thomas Pock
acknowledges support from the Austrian science fund under the
ANR-FWF project \enquote{Efficient algorithms for nonsmooth optimization in
imaging}, No. I1148 and the FWF-START project \enquote{Bilevel optimization
for Computer Vision}, No. Y729.


\bibliographystyle{spmpsci}
\bibliography{ochs}

\end{document}